\documentclass[12pt]{amsart}

\usepackage{amssymb, amsmath}
\usepackage{mathrsfs, bbm}
\usepackage{graphicx, enumerate}
\usepackage[margin=2.9cm]{geometry}
\usepackage{caption}
\usepackage{amsthm}
\usepackage{tikz}
\usepackage{tikz-cd}
\usepackage{svg}

\theoremstyle{definition}
\newtheorem{theorem}{{Theorem}}
\newtheorem{lemma}{{Lemma}}
\newtheorem{proposition}{{Proposition}}
\newtheorem{definition}{{Definition}}
\newtheorem{corollary}{{Corollary}}

\newtheorem{conjecture}{{Conjecture}}
\newtheorem{remark}{{Remark}}
\newtheorem{claim}{{Claim}}
\newtheorem{question}{{Question}}

\newcommand{\bdefn}{\begin{definition}}
\newcommand{\edefn}{\end{definition}}
\newcommand{\bthm}{\begin{theorem}}
\newcommand{\ethm}{\end{theorem}}
\newcommand{\bprop}{\begin{proposition}}
\newcommand{\eprop}{\end{proposition}}
\newcommand{\bcor}{\begin{corollary}}
\newcommand{\ecor}{\end{corollary}}
\newcommand{\blem}{\begin{lemma}}
\newcommand{\elem}{\end{lemma}}
\newcommand{\brmk}{\begin{remark}}
\newcommand{\ermk}{\end{remark}}
\newcommand{\bpf}{\begin{proof}}
\newcommand{\epf}{\end{proof}}
\newcommand{\bclm}{\begin{claim}}
\newcommand{\eclm}{\end{claim}}
\newcommand{\bquest}{\begin{question}}
\newcommand{\equest}{\end{question}}

\newcommand{\bi}{\begin{itemize}}
\newcommand{\ei}{\end{itemize}}
\newcommand{\be}{\begin{enumerate}}
\newcommand{\ee}{\end{enumerate}}

\makeatletter
\newtheorem*{rep@theorem}{\rep@title}
\newcommand{\newreptheorem}[2]{%
\newenvironment{rep#1}[1]{%
 \def\rep@title{#2 \ref{##1}}%
 \begin{rep@theorem}}%
 {\end{rep@theorem}}}
 \newcommand{\newinfreptheorem}[2]{%
\newenvironment{infrep#1}[1]{%
 \def\rep@title{#2 \ref{##1} (informal)}%
 \begin{rep@theorem}}%
 {\end{rep@theorem}}}
\makeatother
\newreptheorem{theorem}{Theorem}
\newreptheorem{corollary}{Corollary}
\newreptheorem{conjecture}{Conjecture}

\newcommand{\map}[1]{\xrightarrow{#1}}
\newcommand{\inv}{^{-1}}
\newcommand{\st}{\ | \ }

\newcommand{\R}{\mathbb{R}}

\newcommand{\C}{\mathbb{C}}

\newcommand{\Z}{\mathbb{Z}}

\newcommand{\CB}{\mathcal{B}}
\newcommand{\CC}{\mathcal{C}}

\newcommand{\CH}{\mathcal{H}}
\newcommand{\CD}{\mathcal{D}} 
\newcommand{\CT}{\mathcal{T}}

\newcommand{\into}{\hookrightarrow}

\newcommand{\bpmat}{\begin{pmatrix}}
\newcommand{\epmat}{\end{pmatrix}}

\newcommand{\Ext}{\mathrm{Ext}}
\newcommand{\F}{\mathbb{F}}

\newcommand{\sltwo}{\mathfrak{sl}_2}

\newcommand{\fix}{^\mathrm{fix}}
\newcommand{\vhE}{\,^{vh}E}
\newcommand{\hvE}{\,^{hv}E}
\newcommand{\hv}{^{hv}}
\newcommand{\rk}{\text{rk}_{\mathbb{F}_2}}
\newcommand{\abuts}{\rightrightarrows}
\newcommand{\atate}{\text{\sc{AKh}}^{\text{Tate}}}
\newcommand{\khtate}{\text{\sc{Kh}}^{\text{Tate}}}
\newcommand{\FF}{\mathbb{F}_2}

\newcommand{\scube}{\text{\sc Cube}} 
\newcommand{\dtwo}{^{hv}d^2}

\newcommand{\SA}{\mathscr{A}}
\newcommand{\SB}{\mathscr{B}}
\newcommand{\SC}{\mathscr{C}}
\newcommand{\SD}{\mathscr{D}}
\newcommand{\SE}{\mathscr{E}}
\newcommand{\SF}{\mathscr{F}}

\newcommand{\bcent}{\begin{center}}
\newcommand{\ecent}{\end{center}}
\newcommand{\bleft}{\begin{flushleft}}
\newcommand{\eleft}{\end{flushleft}}

\newcommand{\bfig}{\begin{figure}\centering
\begin{minipage}{\textwidth}\centering}
\newcommand{\efig}{\end{minipage}
\end{figure}}

\title[Annular Khovanov homology of 2-periodic links]{A rank inequality for the annular Khovanov homology of 2-periodic links}

\author{Melissa Zhang}

\begin{document}

\maketitle

\begin{abstract}
For a 2-periodic link $\tilde L$ in the thickened annulus and its quotient link $L$, we exhibit a spectral sequence with 
$E^1 \cong AKh(\tilde L) \otimes_{\FF} \FF[\theta, \theta\inv] 
\abuts E^\infty \cong AKh(L) \otimes_{\FF} \FF[\theta, \theta\inv].$
This spectral sequence splits along quantum and $\sltwo$ weight space  gradings, proving a rank inequality 
$ \rk AKh^{j,k}(L) \leq \rk  AKh^{2j-k,k} (\tilde L)  $ for every pair of quantum and $\sltwo$ weight space  gradings $(j,k)$. 
We also present a few decategorified consequences and discuss partial results toward a similar statement for the Khovanov homology of 2-periodic links.
\end{abstract}

\tableofcontents

\section{Introduction} \label{introduction}

The world of topological objects is rich with spaces exhibiting interesting symmetries.
These symmetries manifest themselves in the algebraic invariants we use to study the objects and can be manipulated to provide insights into the object itself. 
As fundamental, ubiquitous objects in low-dimensional topology, knots and links are some of the most natural objects to have algebraic invariants associated to them.

The rise of homology-type invariants categorifying polynomial link invariants have led to the study of equivariant homology theories by Chbili \cite{Chbili-2007}, Seidel and Smith \cite{Seidel-Smith-2010}, Hendricks \cite{Hendricks-2012} \cite{Hendricks-2015}, Politarczyk \cite{Politarczyk-2015} \cite{Politarczyk-2017}, Borodzik and Politarczyk \cite{Boro-Polit-2017}, and many others, building and improving upon earlier work of Murasugi \cite{Murasugi-1971} \cite{Murasugi-1988}, Yokota \cite{Yokota-1993}, and Przytycki \cite{Prz-2004} relating the polynomial invariants of periodic links with those of their quotient links.

A link $\tilde L$ in $S^3$ is \emph{$p$-periodic} if there is a $p$-periodic automorphism $\tau$ of $S^3$ with fixed set an unknot $B$ restricting to a diffeomorphism on $\tilde L$; the \emph{quotient link} $L$ is the quotient of $\tilde L$ under the action of $\tau$. 
Since $\tilde L$ can naturally be viewed in a thickened annulus $A\times I$ by deleting the fixed set $B$ from $S^3$, Asaeda, Przytycki, and Sikora's \emph{annular Khovanov homology} \cite{APS-2004} seems particularly suitable for probing periodic links. 
We'll refer to a link in the thickened annulus as an \emph{annular link}.

In this paper, we construct a Tate-like bicomplex (\S \ref{annular-khovanov-tate-bicomplex}) to show rank inequalities for the annular Khovanov homology of a 2-periodic annular link $\tilde L$ and that of its quotient link $L$. For any annular link $L$, let $AKh^{j,k}(L)$ denote the annular Khovanov homology of $L$ at quantum grading level $j$ and $\sltwo$ weight space grading level $k$.
\bthm \label{main-thm}
Let $\tilde L$ be a 2-periodic link with quotient link $L$. For each pair of integers $(j,k)$, there is a spectral sequence with
$$
E^1 \cong AKh^{2j-k,k}(\tilde L) \otimes_{\FF} \FF[\theta, \theta\inv] 
\abuts
E^\infty \cong AKh^{j,k}(L) \otimes_{\FF} \FF[\theta, \theta\inv].
$$
\ethm

A rank inequality immediately follows:

\bcor \label{main-cor}
For $\tilde L$ and $L$ as above, 
$$
\rk AKh^{j,k}(L) \leq \rk AKh^{2j-k, k}(\tilde L).
$$
\ecor

This result generalizes \cite{Cornish-2016}, where Cornish proves a rank inequality between the next-to-top $\sltwo$ weight space gradings in the annular Khovanov homologies of a 2-periodic link and that of its quotient by employing the Lipshitz-Treumann spectral sequence \cite{Lipshitz-Treumann-2016} with Auroux, Grigsby, and Wehrli's identification of the next-to-top winding grading with the Hochschild homology of a Khovanov-Seidel bimodule \cite{Auroux-Grigsby-Wehrli-2015}. He remarks that by using Beliakova-Putyra-Wehrli's generalization of Auroux, Grigsby, and Wehrli's work \cite{BPW-2016}, one should be able to prove similar rank inequalities in other $\sltwo$ weight space gradings, but that this requires checking the Lipshitz-Treumann algebraic conditions for larger dg algebras.

While our result does not require the links to be braid closures, annular Khovanov homology is often used to study braids, as the set of conjugacy classes of a braid group $B_n$ embed into the set of isotopy classes of annular links. In particular, $B_n$ can be viewed as the mapping class group of a disk with $n$ punctures, and Corollary \ref{main-cor} suggests that the rank of annular Khovanov homology could provide a measure of complexity in the mapping class group.

\paragraph{} 
Following the notation in \cite{Roberts-2013} (page 418, under Definition 2.2; also see \S \ref{decat} in this paper), the decategorification of Theorem \ref{main-thm} can be written as follows:
\bcor \label{decat-cor}
For all $j$ and $k$, 
$$\langle q_{k, \tilde L}(-1,q), q^{2j-k} \rangle \equiv \langle q_{k,L}(-1,q), q^j \rangle \mod 2$$
where $\langle f,g \rangle$ denotes the coefficient of $g$ in $f$. 
\ecor

However, the spectral sequence hints that the grading $j_1 := j-k$ is actually more pertinent to annular links (see \cite{GLW-2015} and \cite{GLW-2017} for more evidence). With this in mind, one can write down the decategorification as a statement similar to Murasugi's theorem on the Jones polynomial of periodic links \cite{Murasugi-1988}:
\bcor \label{murasugi-like}
$V_{\tilde L}(1,q, q\inv) \equiv [V_L(1,q,q\inv)]^2 \mod 2.$
\ecor

The Jones polynomial of a periodic link $\tilde L$ is related to the decategorification of the annular Khovanov homology of its quotient link in the following way: 
\bcor \label{decat-kh}
$V_{\tilde L}(1,q,1) \equiv V_L(1,q^2,q\inv) \mod 2.$
\ecor

Moreover, we conjecture that a similar spectral sequence also relates $Kh(\tilde L)$ and $AKh(L)$ (see \S \ref{section-conjecture}):

\begin{conjecture} \label{main-conj}
Let $\tilde L$ be a 2-periodic link in $S^3$ with quotient link $L$. There is a spectral sequence with
$$
E^1 \cong Kh(\tilde L) \otimes_{\FF} \FF[\theta, \theta\inv]
\abuts
E^\infty \cong AKh(L) \otimes_{\FF} \FF[\theta, \theta\inv].
$$
This would in turn imply the following cascade of rank inequalities:
$$
\rk AKh(\tilde L) 
\geq \rk Kh(\tilde L)
\geq \rk AKh(L)
\geq \rk Kh(L)
$$
where the first and third inequalities are given by the $k$-grading filtration on $CKh(\CD(\tilde L))$ and $CKh(\CD(L))$.
\end{conjecture}

In fact, Corollary \ref{decat-kh} would be the decategorification of Conjecture \ref{main-conj}. We can show that the conjecture and consequent rank inequalities hold for the following family of annular links:

\bthm \label{mostly-negative-thm}
If the annular braid closure $L = \widehat \beta$ has a diagram with at most 1 positive crossing, then the spectral sequence in Conjecture \ref{main-conj} exists and the cascade of rank inequalities holds.
\ethm

Furthermore, it follows from Theorem \ref{mostly-negative-thm} that the cascade of rank inequalities also holds for positive braid closures (and more):

\bcor \label{mostly-positive-cor}
If the annular braid closure $L = \widehat \beta$ has a diagram with at most 1 negative crossing, then the cascade of rank inequalities holds. 
\ecor

\subsection{Acknowledgements} \label{acknowledgements}

I am grateful to John Baldwin for suggesting the problem and guiding me; to Eli Grigsby for her constant support, advice, and guidance; to Spencer Leslie and David Treumann for enlightening conversations about equivariant homology and Smith theory; to Robert Lipshitz and Kristen Hendricks for useful discussions about the main result; and to Adam Saltz, Eli, and John for reading and commenting on drafts.

\section{Algebraic preliminaries} \label{algebraic-preliminaries}

We begin by discussing the necessary homological algebra. The main theorem follows from computing two related spectral sequences arising from a particular bicomplex (constructed in \S \ref{annular-khovanov-tate-bicomplex}). Spectral sequences from bicomplexes are a special case of spectral sequences from filtered complexes, which we discuss first. 

The spectral sequences will be computed explicitly by making use of a well-known \emph{cancellation lemma} (Lemma \ref{cancellation-lemma}). Due to a special property of our complexes (see Definition \ref{filtered-graded-basis}), the algebraic computations can be described visually using dots and arrows; we describe this at the end of \S \ref{filtered-chain-complexes}. 

In \S \ref{equiv-hom-tate-ss}, we give a brief overview of the ideas in Borel equivariant cohomology which motivate the constructions in this paper.

\subsection{Filtered chain complexes} \label{filtered-chain-complexes}

Let $\CC_\bullet$ be a chain complex of $\F_2$ vector spaces, with differential $\partial$.

A \emph{decreasing $\Z$-filtration} of $\CC_\bullet$ is a sequence of subcomplexes indexed by decreasing integers:
$$ 
\ldots  \subset
F_n \CC_\bullet \subset F_{n-1} \CC_\bullet \subset F_{n-2} \CC_\bullet \subset \ldots \subset F_m \CC_\bullet \subset \ldots \CC_\bullet.
$$

The filtration $F_* \CC_\bullet$ provides a \emph{filtration grading} for the chains in $\CC_\bullet$: $x \in \CC_\bullet$ has filtration grading $p$ if and only if $x \in F_p \CC_\bullet$ and $x \not\in F_{p+1} \CC_\bullet$. 
The \emph{associated graded vector space} $G_* \CC_\bullet$ is given by
$$ G_p \CC_\bullet = F_p \CC_\bullet / F_{p+1} \CC_\bullet.$$

The existence of the filtration on $\CC_\bullet$ implies that every component of $\partial$ originating from $x \in \CC_\bullet$ maps to targets $y$ where the filtration grading of $y$ is at least that of $x$.
For each $s \in \Z$, define $\partial_s$ to be the sum of the components of $\partial$ which \emph{shift the filtration grading} by $s$, or are of  \emph{degree $s$ with respect to the filtration grading}. 
In other words, $\partial_s = 0 $ for all $s < 0$, and we say that $\partial$ is non-negative with respect to the filtration grading.

Let $\CH_\bullet$ denote the homology of the chain complex $\CC_\bullet$.
The filtration $F_*\CC_\bullet$ induces a filtration $F_* \CH_\bullet$ on $\CH_\bullet$:
the class $[x] \in \CH_\bullet$ is in $F_p \CH_\bullet$ if and only if it is represented by some chain $x \in [x]$ where $x \in F_p \CC_\bullet$.
The induced filtration gives rise to the associated graded vector space for homology, $G_* \CH_\bullet$.

\bdefn \label{filtered-graded-basis}

In our work, the underlying vector space of $\CC_\bullet$ is always freely generated by a finite collection of distinguished generators which come equipped with some $\Z$-grading $gr$, and the differential is non-negative with respect to $gr$.
In this case, $gr$ \emph{induces} a $\Z$-filtration and takes on the role of the filtration grading. We call the distinguished generators a \emph{filtered graded basis}. 

Similarly, a bigraded underlying vector space induces a $\Z \oplus \Z$-filtration if the differential is monotone with respect to each grading.

\edefn

\begin{figure} 
\centering 
\def\svgwidth{200pt} 
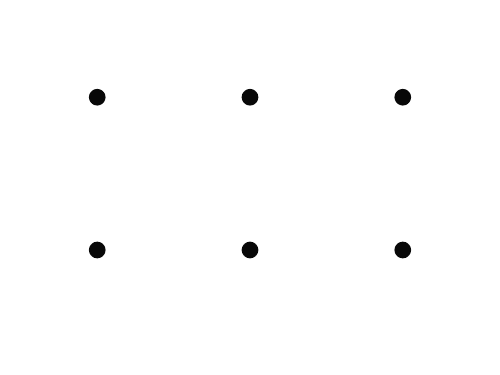
\caption{}
\label{dots-n-arrows} 
\end{figure}

Let $(\CC, \partial)$ be such a complex. 
We can visualize $(\CC, \partial)$ by using dots to represent the distinguished generators and arrows to represent components of the differential, as shown in Figure \ref{dots-n-arrows}. 
We write $\partial(x,y)$ for the coefficient of $y$ in $\partial(x)$, so that $\partial(x)$ is the sum of $\partial(x,y)\cdot y$ as $y$ runs over all the distinguished generators of $\CC$. 
If $\partial(x,y) = 1$, we draw an arrow from $x$ to $y$.

\subsection{Spectral sequences from filtered complexes} \label{spectral-sequences-from-filtered-complexes}

\begin{figure} 
\centering 
\def\svgwidth{300pt} 
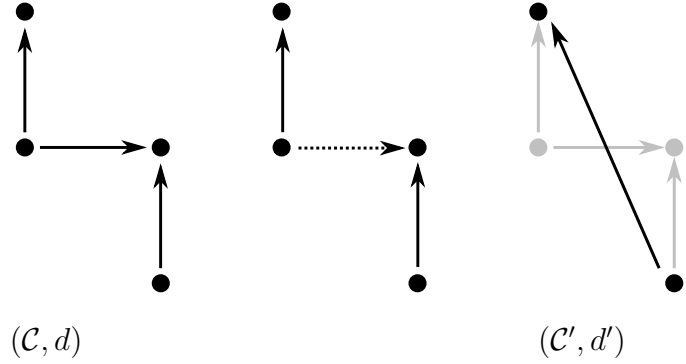
\caption{The cancellation process. On the left is the original complex. In the middle, the dotted arrow represents the arrow about to be cancelled. On the right is the resulting complex, with fewer generators and arrows, but with the same homology as the original complex.}
\label{cancellation} 
\end{figure}

Spectral sequences are used to gradually approximate the true homology of a complex by gradually simplifying the complex, whilst preserving the homology. 
We execute this by repeatedly applying the \emph{cancellation} operation  (Figure \ref{cancellation}), which relies on the following well-known ``cancellation lemma.''

\blem \label{cancellation-lemma}
[\cite{Baldwin-2011}, Lemma 4.1]

Let $(\CC, \partial)$ be a complex of $\F_2$ vector spaces freely generated by elements $x_i$. Let $\partial(x_i,x_j)$ be the coefficient of $x_j$ in $\partial(x_i)$. If $\partial(x_k, x_l) = 1$, define a new complex $(\CC', \partial')$ with generators $\{x_i \st i \neq k,l\}$ and differential
$$
\partial'(x_i) = \partial(x_i) + \partial(x_i,x_l) \partial(x_k).$$
Then $(\CC', \partial')$ is chain homotopy equivalent to $(\CC, \partial)$.
\elem

Given a decreasing $\Z$-filtered complex $F_*\CC_\bullet$ with filtered differential $\partial$, we can build a spectral sequence to compute the homology by starting with the associated graded complex and sequentially cancelling components of the differential of larger and larger filtration degree.

The pages of our spectral sequence are denoted $(E^r, d^r)$, where $E^r_p$ is the vector space at filtration grading $p$ surviving to page $r$, $E^r = \oplus_p E^r_p$, and $d^r$ is the induced differential on $E^r$ induced by cancellation on the previous pages. For example, here are the first few pages:

\bi
\item On the 0th page, $E^0_p = G_p\CC / G_{p-1} \CC$; $d^0$ is the sum of all the components of $\partial$ which preserve the filtration grading $gr$.

\item Canceling all the components of $d^0$, we obtain a chain homotopy equivalent complex $(\CC', \partial')$ where the underlying vector space $\CC'$ injects into $\CC$ and therefore has an induced filtration. Since there are no longer differentials which preserve the filtration grading, all components of $\partial'$ now shift the filtration grading by at least 1. We call $E^1_p := G_p\CC' / G_{p-1}\CC'$ the vector space at filtration grading $p$ that \emph{survives to page 1}, and $d^1$ is the sum of the components of $\partial'$ which shift the filtration grading by exactly 1.

\item Canceling all the components of $d^1$ in the complex $(\CC', \partial')$, we obtain a chain homotopy equivalent complex $(\CC'', \partial'')$. We still have the filtration grading, since $\CC'' \into \CC' \into \CC$, and we let $d^2$ be the sum of the components of $\partial''$ which shift the grading by exactly 2.
\ei

Iterating this process, we obtain pages $(E^r, d^r)$ for all $r \geq 0$. In some situations, the spectral sequence eventually \emph{collapses}, i.e.\ there is some $N$ such that for all $r \geq N$, $d^r=0$, and all pages are identical from then on. In such cases, the \emph{limit term} $E^\infty := E^N = E^{N+1} = \ldots$ is identified with the homology of the original complex $(\CC, \partial)$. We say that the spectral sequence $E^\bullet$ \emph{abuts to} $E^\infty$. 

\begin{figure} 
\centering 
\def\svgwidth{400pt} 
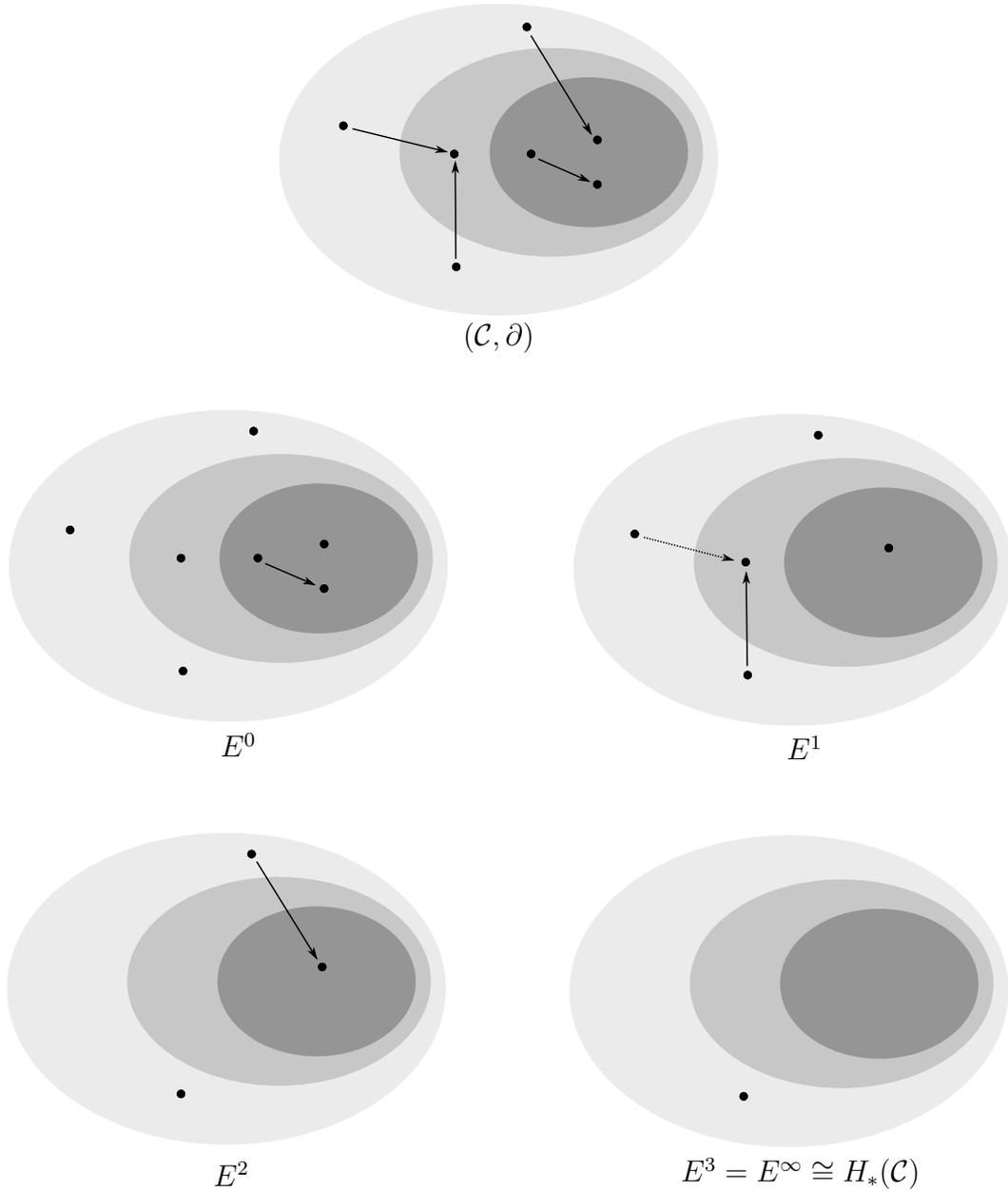
\caption{The spectral sequence which computes the homology of a filtered complex.}
\label{filtration-ss-example} 
\end{figure}

Figure \ref{filtration-ss-example} shows the cancellation process (the spectral sequence) used to compute the homology of filtered complex.

\subsection{Bicomplexes} \label{bicomplexes}

\bfig
\begin{tikzcd}
 & \vdots & \vdots & \vdots &  \\
\cdots  \arrow{r}{\partial^h} 
	& \CC_{p-1,q+1} \arrow{u}{\partial^v} \arrow{r}{\partial^h} 
	& \CC_{p,q+1} \arrow{u}{\partial^v} \arrow{r}{\partial^h} 
	& \CC_{p+1,q+1} \arrow{u}{\partial^v} \arrow{r}{\partial^h} 
	& \cdots \\
\cdots  \arrow{r}{\partial^h} 
	& \CC_{p-1,q} \arrow{u}{\partial^v} \arrow{r}{\partial^h} 
	& \CC_{p,q} \arrow{u}{\partial^v} \arrow{r}{\partial^h} 
	& \CC_{p+1,q} \arrow{u}{\partial^v} \arrow{r}{\partial^h} 
	& \cdots \\
\cdots  \arrow{r}{\partial^h} 
	& \CC_{p-1,q-1} \arrow{u}{\partial^v} \arrow{r}{\partial^h} 
	& \CC_{p,q-1} \arrow{u}{\partial^v} \arrow{r}{\partial^h} 
	& \CC_{p+1,q-1} \arrow{u}{\partial^v} \arrow{r}{\partial^h} 
	& \cdots \\
	& \vdots \arrow{u}{\partial^v} 
	& \vdots \arrow{u}{\partial^v}
	& \vdots \arrow{u}{\partial^v}
	& %
\end{tikzcd}
\caption{}
\label{bicomplex}
\efig

A \emph{bicomplex} $\CC_{\bullet, \bullet}$ (Figure \ref{bicomplex}) is a bigraded $\F_2$ vector space with differentials
$$ \partial_{p,q}^h: \CC_{p,q} \to \CC_{p+1,q} $$
$$  \partial_{p,q}^v: \CC_{p,q} \to \CC_{p,q+1}$$
such that $\partial^h \circ \partial^v + \partial^v \circ \partial^h = 0$.
Here $\partial^h$ is the sum of all of the \emph{horizontal differentials} $\partial_{p,q}^h$ and $\partial^v$ is the sum of all of the \emph{vertical differentials} $\partial_{p,q}^v$.

The corresponding \emph{total complex} $Tot(\CC)_\bullet$ is given by
$$ Tot(\CC)_n = \bigoplus_{p+q = n} \CC_{p,q}$$
with differential $\partial^{Tot} = \partial^h + \partial^v$.

\bfig
\begin{minipage}{0.47\textwidth}
\def\svgwidth{\textwidth} 
\begingroup%
  \makeatletter%
  \providecommand\color[2][]{%
    \errmessage{(Inkscape) Color is used for the text in Inkscape, but the package 'color.sty' is not loaded}%
    \renewcommand\color[2][]{}%
  }%
  \providecommand\transparent[1]{%
    \errmessage{(Inkscape) Transparency is used (non-zero) for the text in Inkscape, but the package 'transparent.sty' is not loaded}%
    \renewcommand\transparent[1]{}%
  }%
  \providecommand\rotatebox[2]{#2}%
  \ifx\svgwidth\undefined%
    \setlength{\unitlength}{371.34723333bp}%
    \ifx\svgscale\undefined%
      \relax%
    \else%
      \setlength{\unitlength}{\unitlength * \real{\svgscale}}%
    \fi%
  \else%
    \setlength{\unitlength}{\svgwidth}%
  \fi%
  \global\let\svgwidth\undefined%
  \global\let\svgscale\undefined%
  \makeatother%
  \begin{picture}(1,0.77770873)%
    \put(0,0){\includegraphics[width=\unitlength,page=1]{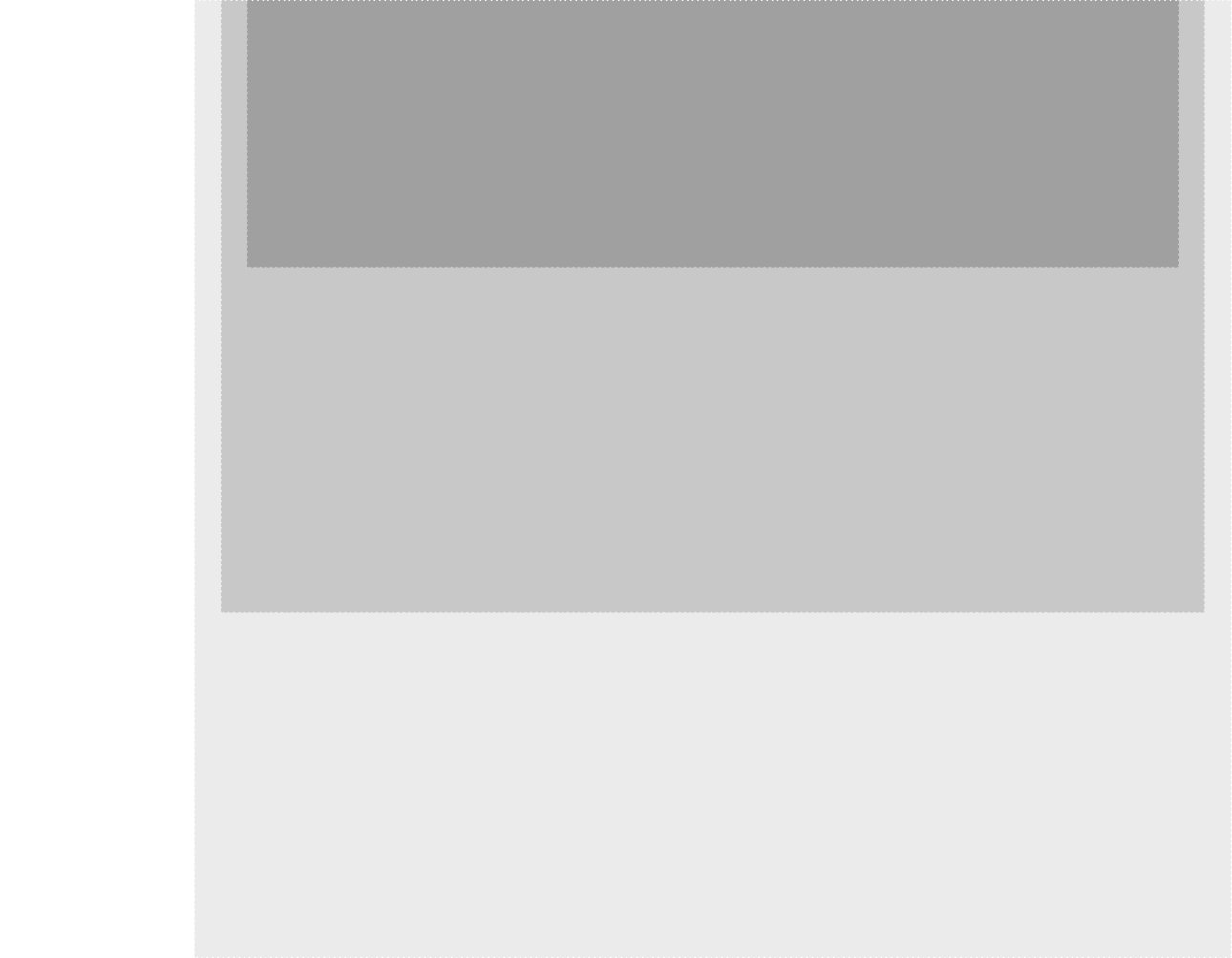}}%
    \put(-0.00163047,0.70561167){\color[rgb]{0,0,0}\makebox(0,0)[lb]{\smash{$^v F_2$}}}%
    \put(-0.00163047,0.42555033){\color[rgb]{0,0,0}\makebox(0,0)[lb]{\smash{$^v F_1$}}}%
    \put(-0.00163047,0.12394585){\color[rgb]{0,0,0}\makebox(0,0)[lb]{\smash{$^v F_0$}}}%
    \put(0.29237924,0.61960207){\color[rgb]{0,0,0}\makebox(0,0)[lb]{\smash{$\mathcal{C}_{0, 2}$}}}%
    \put(0.19105371,0.69045891){\color[rgb]{0,0,0}\makebox(0,0)[lt]{\begin{minipage}{0.44163516\unitlength}\raggedright \end{minipage}}}%
    \put(0.29237924,0.33954074){\color[rgb]{0,0,0}\makebox(0,0)[lb]{\smash{$\mathcal{C}_{0, 1}$}}}%
    \put(0.29237924,0.05947934){\color[rgb]{0,0,0}\makebox(0,0)[lb]{\smash{$\mathcal{C}_{0, 0}$}}}%
    \put(0.55089741,0.61960207){\color[rgb]{0,0,0}\makebox(0,0)[lb]{\smash{$\mathcal{C}_{1, 2}$}}}%
    \put(0.55089741,0.33954074){\color[rgb]{0,0,0}\makebox(0,0)[lb]{\smash{$\mathcal{C}_{1, 1}$}}}%
    \put(0.55089741,0.05947934){\color[rgb]{0,0,0}\makebox(0,0)[lb]{\smash{$\mathcal{C}_{1, 2}$}}}%
    \put(0.80941556,0.61960207){\color[rgb]{0,0,0}\makebox(0,0)[lb]{\smash{$\mathcal{C}_{2, 2}$}}}%
    \put(0.80941556,0.33954074){\color[rgb]{0,0,0}\makebox(0,0)[lb]{\smash{$\mathcal{C}_{2, 1}$}}}%
    \put(0.80941556,0.05947934){\color[rgb]{0,0,0}\makebox(0,0)[lb]{\smash{$\mathcal{C}_{2, 0}$}}}%
  \end{picture}%
\endgroup%

\end{minipage}
\hspace{.04\textwidth}
\begin{minipage}{0.47\textwidth}
\def\svgwidth{.8\textwidth} 
\begingroup%
  \makeatletter%
  \providecommand\color[2][]{%
    \errmessage{(Inkscape) Color is used for the text in Inkscape, but the package 'color.sty' is not loaded}%
    \renewcommand\color[2][]{}%
  }%
  \providecommand\transparent[1]{%
    \errmessage{(Inkscape) Transparency is used (non-zero) for the text in Inkscape, but the package 'transparent.sty' is not loaded}%
    \renewcommand\transparent[1]{}%
  }%
  \providecommand\rotatebox[2]{#2}%
  \ifx\svgwidth\undefined%
    \setlength{\unitlength}{265.44322577bp}%
    \ifx\svgscale\undefined%
      \relax%
    \else%
      \setlength{\unitlength}{\unitlength * \real{\svgscale}}%
    \fi%
  \else%
    \setlength{\unitlength}{\svgwidth}%
  \fi%
  \global\let\svgwidth\undefined%
  \global\let\svgscale\undefined%
  \makeatother%
  \begin{picture}(1,1.18221217)%
    \put(0,0){\includegraphics[width=\unitlength,page=1]{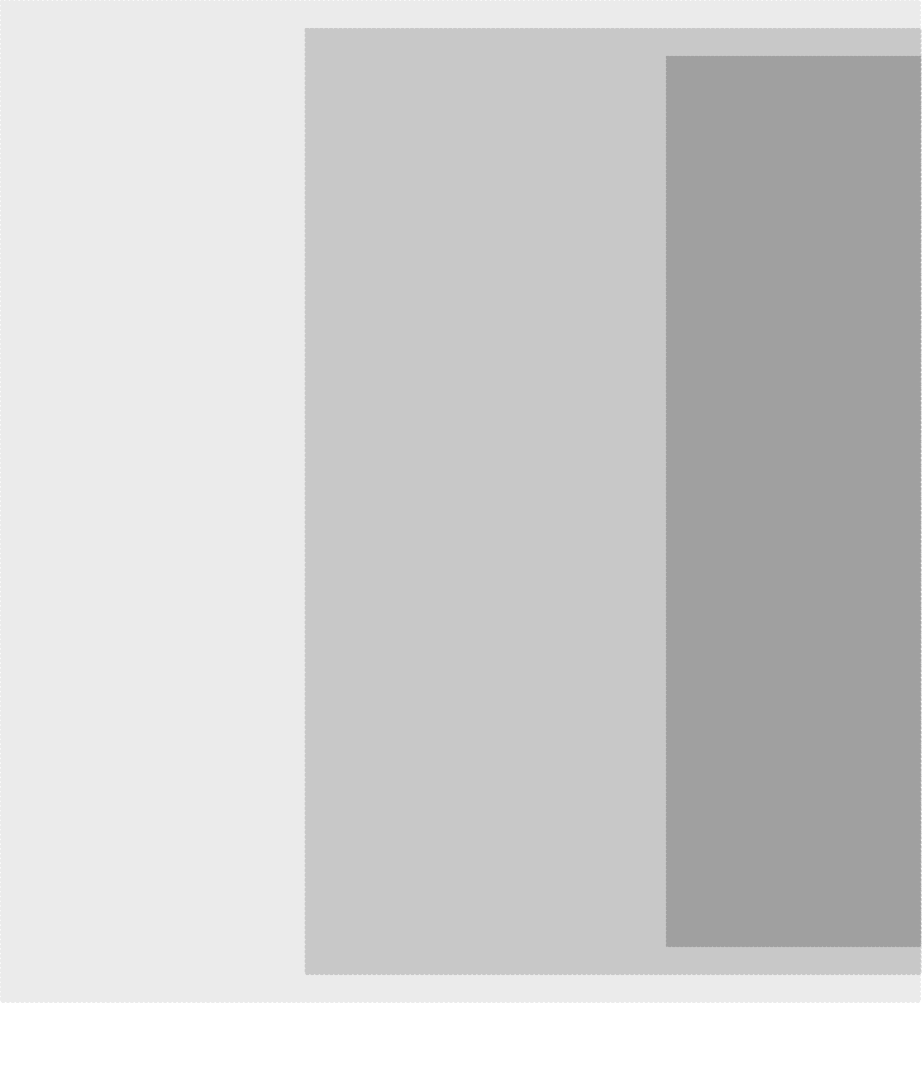}}%
    \put(0.060864,0.00531235){\color[rgb]{0,0,0}\makebox(0,0)[lb]{\smash{$^h F_0$}}}%
    \put(0.42252326,0.00531235){\color[rgb]{0,0,0}\makebox(0,0)[lb]{\smash{$^h F_1$}}}%
    \put(0.81432079,0.00531235){\color[rgb]{0,0,0}\makebox(0,0)[lb]{\smash{$^h F_2$}}}%
    \put(0.06699322,1.02130208){\color[rgb]{0,0,0}\makebox(0,0)[lb]{\smash{$\mathcal{C}_{0, 2}$}}}%
    \put(-0.0747582,1.1204287){\color[rgb]{0,0,0}\makebox(0,0)[lt]{\begin{minipage}{0.61783454\unitlength}\raggedright \end{minipage}}}%
    \put(0.06699322,0.62950455){\color[rgb]{0,0,0}\makebox(0,0)[lb]{\smash{$\mathcal{C}_{0, 1}$}}}%
    \put(0.06699322,0.23770693){\color[rgb]{0,0,0}\makebox(0,0)[lb]{\smash{$\mathcal{C}_{0, 0}$}}}%
    \put(0.4286525,1.02130208){\color[rgb]{0,0,0}\makebox(0,0)[lb]{\smash{$\mathcal{C}_{1, 2}$}}}%
    \put(0.4286525,0.62950455){\color[rgb]{0,0,0}\makebox(0,0)[lb]{\smash{$\mathcal{C}_{1, 1}$}}}%
    \put(0.4286525,0.23770693){\color[rgb]{0,0,0}\makebox(0,0)[lb]{\smash{$\mathcal{C}_{1, 2}$}}}%
    \put(0.79031175,1.02130208){\color[rgb]{0,0,0}\makebox(0,0)[lb]{\smash{$\mathcal{C}_{2, 2}$}}}%
    \put(0.79031175,0.62950455){\color[rgb]{0,0,0}\makebox(0,0)[lb]{\smash{$\mathcal{C}_{2, 1}$}}}%
    \put(0.79031175,0.23770693){\color[rgb]{0,0,0}\makebox(0,0)[lb]{\smash{$\mathcal{C}_{2, 0}$}}}%
  \end{picture}%
\endgroup%

\end{minipage}
\caption{}
\label{two-std-filtrations}
\efig

The total complex has two standard filtrations $^vF$ (the \emph{row-wise filtration}) and $^hF$ (the \emph{column-wise filtration}) induced by the rows and columns of the bicomplex, respectively, defined below, and depicted in Figure \ref{two-std-filtrations}:
$$
(^vF_m\CC)_{p,q}
= 
\begin{cases}
\CC_{p,q} & \text{if } q \leq m \\
0  & \text{otherwise} 
\end{cases}
$$
$$
(^hF_m\CC)_{p,q}
= 
\begin{cases}
\CC_{p,q} & \text{if } p \leq m \\
0  & \text{otherwise} 
\end{cases}.
$$

Indeed, with respect to both the row- and column-wise filtrations, the differential shifts the filtration degree by a nonnegative number.

\subsection{Spectral sequences from bicomplexes} \label{spectral-sequences-from-bicomplexes}

The standard filtrations $^v F$ and $^h F$ induce spectral sequences of filtered complexes $\hvE$ and $\vhE$, respectively. 

For example, consider the row-wise filtration $^v F$ shown on the left of Figure \ref{two-std-filtrations}. 
The pages of the attached spectral sequence $\hvE$ approximate $H_*(Tot(\CC))$ via the following process.
On page $\hvE^0$, the differential $d^0$ consists of all the components of $\partial^{Tot}$ which preserve the filtration grading; that is, $d^0 = \partial^h$.
The induced differential on the next page $(\partial^{Tot})'$ now only has components which shift the filtration grading by $1$ or more. 

The following proposition says that in certain situations, both spectral sequences eventually collapse to $\hvE^\infty \cong \vhE^\infty \cong  H_*(Tot(\CC))$, where the vector space generated by the surviving basis elements at the level set with filtration level $p+q = m$ is identified with $H_m(Tot(\CC))$. A proof may be found in \cite{McCleary-2001}.

\bprop \label{hvE-and-vhE}
If for each $n$, the number of $p$ such that $\CC_{p,n-p} \neq 0$ is finite, then 
$$
\hvE^\infty_{p,q} = ^v F_q H_{p+q} (Tot(\CC)) / ^v F_{q+1} H_{p+q} (Tot(\CC))
$$
$$
\vhE^\infty_{p,q} = ^h F_p H_{p+q} (Tot(\CC)) / ^h F_{p+1} H_{p+q} (Tot(\CC)).
$$
Both spectral sequences converge to the total homology. 
\eprop

\subsection{Equivariant homology and the Tate spectral sequence} \label{equiv-hom-tate-ss}

Our bicomplex will be modeled on the bicomplex used to compute the Borel equivariant cohomology of a topological space $X$ with an involution $\tau$:
$$ 
H^*_{\Z/2\Z} (X; \F_2) : = \Ext_{\F_2[\Z/2\Z]} (C_*(X), \F_2).
$$
Here $C_*(X)$ is the singular chain complex for $X$ with $\F_2$ coefficients.
Since we are working over $\F_2$ coefficients, $(1+\tau)^2 = 0$, so this is indeed a complex. Moreover, one can check that the differentials of $C_*(X)$ and the induced involution $\tau^\#$ commute, and therefore so do their duals. Thus we can build the double complex
$$
0 \to C^*(X; \F_2) \map{1+\tau^\#} C^*(X; \F_2) \map{1+\tau^\#} C^*(X; \F_2) \map{1+\tau^\#} \cdots.
$$ 
which computes the Borel equivariant cohomology of $X$. Equipped with the row-wise and column-wise filtration, this is a bicomplex.

We can view the underlying vector space of the bicomplex as a module 
$C^*(X; \F_2) \otimes_{\F_2} \F_2[\theta]$, where $\theta$ shifts the column-wise filtration degree by $1$ (to the right).

Localizing at $\theta$ gives the bi-infinite \emph{Tate bicomplex} 
$$
C^*(X; \F_2) \otimes_{\F_2} \F_2[\theta, \theta\inv]
= \left ( 
 \cdots \map{1+\tau^\#} C^*(X; \F_2) \map{1+\tau^\#}  C^*(X; \F_2) \map{1+\tau^\#} \cdots
\right )
$$ 
whose total homology, under some finiteness conditions, is isomorphic to $H^*(X^{\fix}; \F_2) \otimes_{\F_2} \F_2[\theta, \theta\inv]$. Here $X^{\fix}$ is the $\tau$-invariant topological subspace of $X$.
See Section 2 in \cite{Lipshitz-Treumann-2016} for more details on $\Z/2\Z$-localization in Borel equivariant cohomology.

\section{Topological preliminaries} \label{topological-preliminaries}

The purpose of this section is to give an overview of the \emph{annular Khovanov homology} of a link while establishing the notation to be used in the proof of the main result. 
In \S \ref{annular-links-tangles} we define annular links and relate them to braid and tangle closures. 
As annular Khovanov homology is computed from an annular link diagram, in \S \ref{2-periodic-links-diagrams} we set up the annular link diagram such that the topological involution on the link translates nicely to an involution on the diagram. 
\S \ref{annular-khovanov-homology} reviews the construction of the annular Khovanov chain complex via a \emph{cube of resolutions} and defines the three gradings attached to the complex.

Throughout, we work in the smooth category with oriented links and $\FF$ coefficients.

\subsection{Annular links and tangles} \label{annular-links-tangles}

\begin{figure} 
\centering 
\def\svgwidth{200pt} 
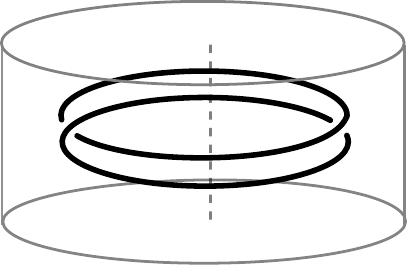
\caption{A Hopf link embedded in a thickened annulus.}
\label{embedded-hopf} 
\end{figure}

An \emph{annular link} is a finite disjoint union of embedded circles $\coprod_1^n S^1$ in a thickened annulus $A \times I$. See Figure \ref{embedded-hopf} for an example.

Annular links are \emph{tangle closures}. A a \emph{tangle} $\CT$ is a properly embedded compact 1-manifold in $D^2 \times I$,
with $\partial \CT \subset D^2 \times \{0,1\}$.
If the number of boundary components of $\CT$ on the two disks agree, then after arranging for $\partial \CT = \{ x_1, \ldots, x_n\} \otimes \{0,1\}$, we can glue $D^2 \times \{0\}$ to $D^2 \times \{1\}$ (via $id_{D^2}$) to obtain a link $\widehat \CT = \CT / (x_i, 0) \sim (x_i,1)$ in a solid torus $D^2 \times I / (D^2 \times \{0\} \sim D^2 \times \{1\})$.

Braids constitute a well-studied case of tangles. Closures of $n$-strand braids form the braid group $B_n$. Isotopy classes of $n$-braid closures correspond to the conjugacy classes in $B_n$. 
With this important case in mind, we think of $A$ as the $xy$-plane punctured at the origin $*$ and call $\{*\} \times I \subset \R^2 \times I \cong D^2 \times I$ the \emph{braid axis}.

\begin{figure} 
\centering 
\def\svgwidth{200pt} 
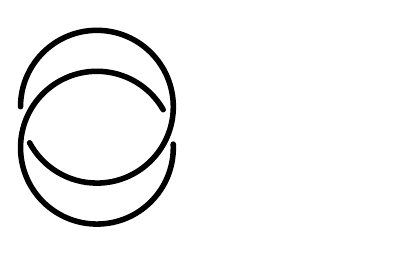
\caption{On the left is an annular diagram for the Hopf link from Figure \ref{embedded-hopf}. On the right is an annular diagram for the quotient knot, a positively stabilized annular unknot.}
\label{annular-hopf} 
\end{figure}

A generic projection of an annular link $L \subset A \times I$ to $A$  is an immersion with at most double points, as in Figure \ref{annular-hopf}. If we also remember the crossing information (which strand was over the other before the projection), then we call the image of this projection a \emph{diagram} of $L$ in $A$, denoted $\CD(L)$. By abuse of notation, sometimes we'll think of $\CD(\cdot)$ as the projection $A \times I \to A$, and refer to the image $\CD(L)\subset A$ as the \emph{annular link diagram}, even though the image doesn't actually capture the crossing information.

\subsection{2-periodic links and their diagrams} \label{2-periodic-links-diagrams}

Given a 2-periodic link $\tilde L \subset S^3$ with involution $\tau: S^3 \to S^3$, the fixed point set $B$ (for ``braid axis'') is an unknotted $S^1$ in $S^3$.

With some choice of coordinates on $S^3$ minus a point $\{\infty\}$ on $B$, we can view
$B$ as the $z$-axis in $\R^3$,
$\tau$ as the rotation of $180^\circ$ about the $z$-axis, and let $A$ be the $xy$-plane minus the origin.
Hence $\tilde L \into A \times \R \cong A \times I$ is an annular link.

\brmk
We define a 2-periodic link $\tilde L \subset S^3$ as one that comes with a chosen involution of $S^3$ which fixes an unknot in $S^3$. It is important to note that a link $\tilde L$ in $S^3$ may come with multiple involutions $\tau_i$, in the sense that the link $\tilde L \cup B_1$ is not isotopic to $\tilde L \cup B_2$ (where $B_i$ is the fixed-point set of involution $\tau_i$). In this case $\tilde L$ with involution $\tau_1$ defines a different annular link from $\tilde L$ with involution $\tau_2$. 
\ermk

Maintaining symmetry under $\tau$, isotope $\tilde L$ so that the projection $\CD: A \times I \to A$ takes $\tilde L$ to a link diagram $\CD(\tilde L)$ in $A$. At the double points, we keep the crossing information.
Employing a small isotopy, we may also assume the double points occur away from the $y$-axis, so that there is a clear notion of the ``right side'' of $A \times I$ or $A$ (where $x \geq 0$), and the ``left side''  of $A \times I$ or $A$ (where $x \leq 0$). 

The preimage of the right side of $\CD(\tilde L)$ is a tangle $\CT$, which by symmetry has the same number of loose strands on both ends.
The quotient link of $\tilde L$ with respect to the involution $\tau$  is the tangle closure $L := \widehat \CT$.
A diagram for $L$ in $A$ is given by $\CD(\CT)$ on the right side of $A$ together with an identity braid on the left side of $A$. See Figure \ref{annular-hopf} for example.

\subsection{Annular Khovanov homology} \label{annular-khovanov-homology}

Annular Khovanov homology (with $\F_2$ coefficients) is a TQFT from properly embedded 1-manifolds in $A \times I$ and cobordisms between them to $\F_2$ vector spaces and linear maps between them.
Asaeda, Przytycki, and Sikora define this annular link invariant in \cite{APS-2004}, where they construct a \emph{cube of resolutions} from the diagram $\CD(L)$ and associate to it a triply-graded chain complex $CKh(\CD(L))$ which computes their annular link invariant $AKh(L)$.

\begin{figure} 
\centering 
\def\svgwidth{300pt} 
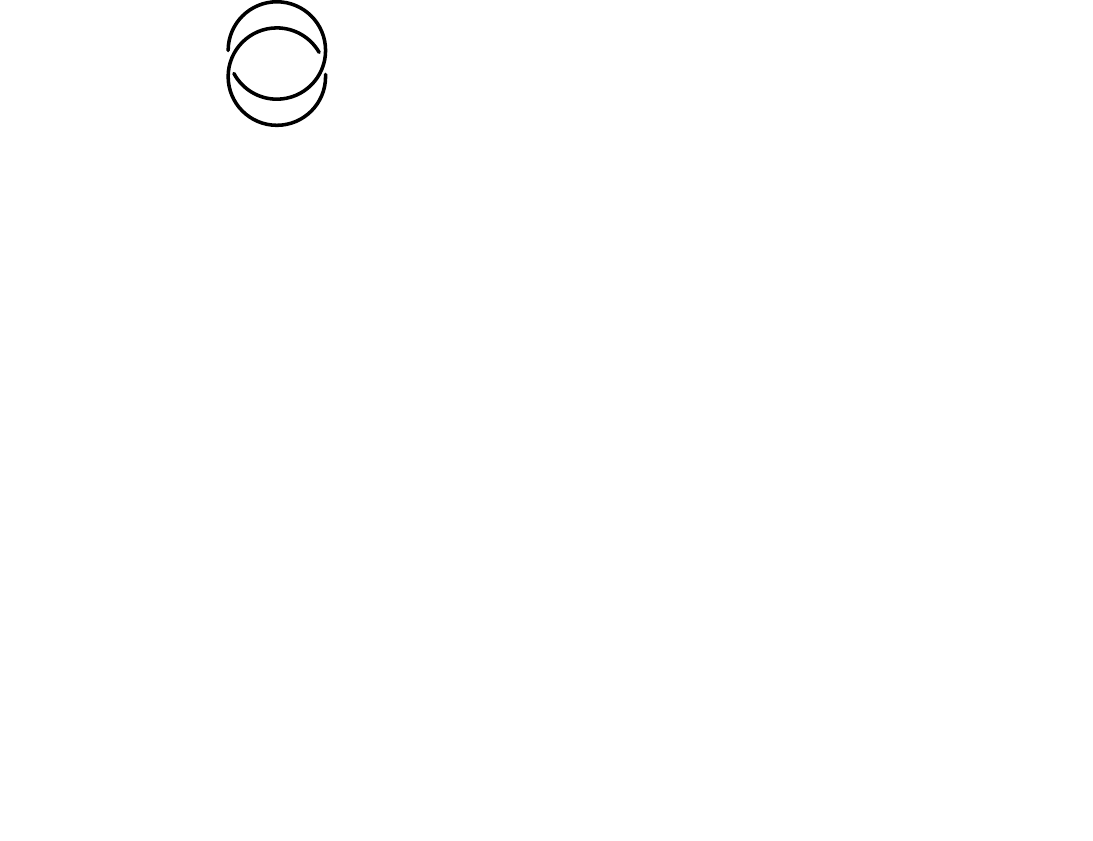
\caption{The cube of resolutions for the annular Hopf link in Figure \ref{annular-hopf}.}
\label{hopf-res-cube} 
\end{figure}

\begin{figure} 
\centering 
\def\svgwidth{400pt} 
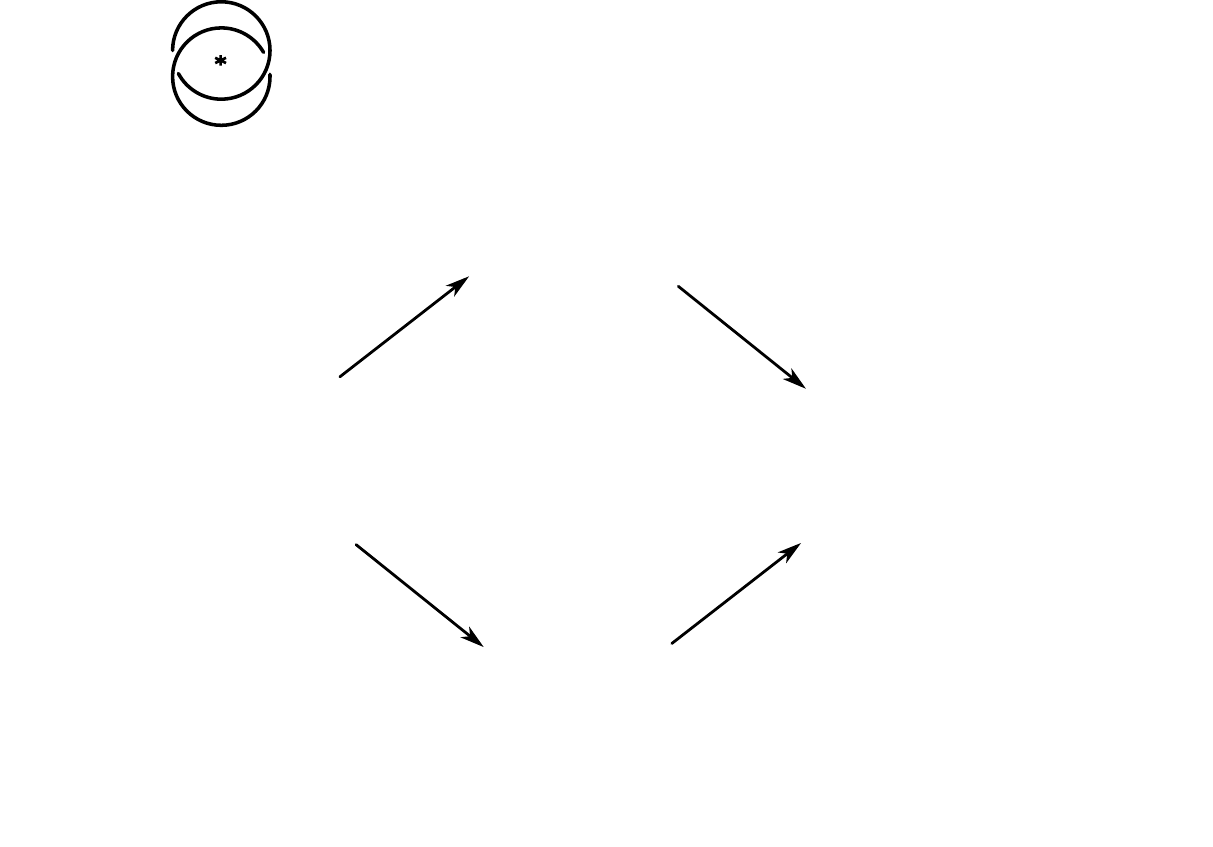
\caption{The cube of chains for the annular Hopf link in Figure \ref{annular-hopf}.}
\label{hopf-vs-cube} 
\end{figure}

We define two  types of cubes, related by the TQFT. See Figure \ref{hopf-res-cube} and Figure  \ref{hopf-vs-cube} for examples.

\bdefn[Cube of resolutions] \label{defn-cube-of-resolutions}
We will use the notation $Cube(\CD)$ to denote the \emph{cube of resolutions} of the diagram $\CD$, whose 
\bi
\item vertices are resolutions of $\CD$, which are closed 1-manifolds embedded in $A$
\item edges are 2-dimensional saddle cobordisms between two vertices.
\ei
\edefn
 
\bdefn[Cube of chains] \label{defn-cube-of-chains}
The annular Khovanov TQFT takes $Cube(\CD)$ to $\text{\sc Cube}(\CD)$, the \emph{cube of chains}, whose 
\bi
\item vertices are $\F_2$ vector spaces
\item edges are linear maps between the vertices.
\ei
\edefn

We outline the construction of $Cube(\CD)$ and $\text{\sc Cube}(\CD)$ in the rest of this subsection.
Let $n$ be the number of crossings in $\CD(L)$, and label the crossings $1$ through $n$ in any order.
Let $n_+$ denote the number of positive crossings, and $n_-$ the number of negative crossings.

\bdefn[Bitstring notation] \label{bitstring-notation}

We call $\alpha \in \{ 0,1 \}^n$ a \emph{bitstring of length $n$}.
There are many operations we can perform on $\alpha$:

\bi
\item $\alpha[i]$ denotes the bit at the $i$-th position, where indexing begins at $i = 0$.
\item $|\alpha| = \sum_{i=0}^{n-1} \alpha[i]$ is the \emph{Hamming weight} of $\alpha$.

\item Let $\alpha' \in \{0,1\}^n$. Suppose there is an index $k$ such that
	\bi
	\item $\alpha[k] = 0$
	\item $\alpha'[k] = 1$
	\item for all $j \neq k$, $\alpha[j] = \alpha'[j]$.
	\ei
In this case, we say that $\alpha'$ is a \emph{bit increment from} $\alpha$. 
\item Let $\beta \in \{0,1\}^n$. Bitstrings of length $n$ form a poset: we write $\alpha < \beta$ if for all $i$, $\alpha[i] \leq \beta[i]$. 
\item $\alpha\beta$ represents the concatenation of $\alpha$ and $\beta$, in that order. That is, $\alpha\beta$ is the length $2n$ bitstring where 
$$
\alpha\beta[i] = 
\begin{cases} 
     \alpha[i] & \text{ for } i < n  \\ \beta[i-n] & \text{ for } i \geq n
\end{cases} .
$$
\item For $\alpha\beta$ as above, there is an involution $\tau$ defined by $\tau(\alpha\beta) = \beta \alpha$.
\ei

\edefn

\brmk \label{tau-remark}
We will also use the symbol ``$\tau$'' to denote some topological involutions, but the distinction from this bitstring operation should be clear from context. See \S \ref{involution-on-ckh} for more on the notation for various involutions.
\ermk

\subsubsection*{Vertices} \label{vertices}

\begin{figure} 
\centering 
\def\svgwidth{300pt} 
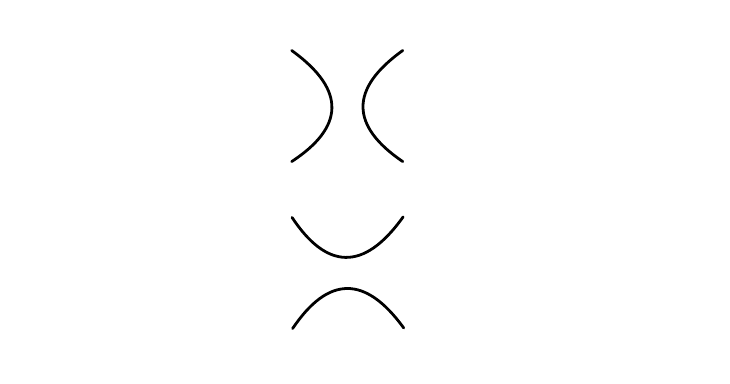
\caption{}
\label{resolution-rule} 
\end{figure}

Let $\alpha \in\{0,1\}^n$ be a length $n$ bitstring.
Figure \ref{resolution-rule} shows two ways of resolving a local picture of a crossing.
Associate to $\alpha$ a resolution $R_\alpha$ of $\CD(L)$, where at crossing $i$ we replace a local picture with a {0-resolution} if $\alpha[i]= 0$, or a {1-resolution} if $\alpha[i] = 1$.
The resulting resolution $R_\alpha$ is a closed 1-manifold embedded in the annulus $A$.

Let $|R_\alpha|$ denote the number of components of $R_\alpha$.
Label the components $C_1, \ldots, C_{|R_\alpha|}$.
Let $\mathbb{V} = \langle v_-, v_+ \rangle$ be the 2-dimensional $\F_2$ vector space generated by the symbols $v_-$ and $v_+$.
We associate to the resolution $R_\alpha$ a vector space $V_\alpha = \mathbb{V}^{\otimes |R_\alpha|}$, generated by the pure tensors of the form $ v = v_{p_1} \otimes v_{p_2} \otimes \ldots \otimes v_{p_{|R_\alpha|}} \in \{ v_-, v_+ \}^{\otimes |R_\alpha|}$. 
We endow $CKh(\CD(L))$ with the following three gradings by assigning gradings to these distinguished generators. This notation is slightly modified from that in \cite{GLW-2017}.

\bi
\item The (homological) \emph{$i$-grading} of $ v$ is 
$$gr_i( v) = |\alpha| - n_-.$$ 

\item Define $p(v) = \# \{ p_i: p_i = +\} -  \# \{ p_i: p_i = -\}$. The (quantum) \emph{$j$-grading} of $ v$ is 
$$ gr_j( v) = p(v)+ |\alpha| + n_+ - 2 n_-.$$

\item The components of $R_\alpha$ come in two flavors with respect to the basepoint $*$:
\bi
\item The \emph{nontrivial circles} have nonzero winding number with respect to $*$ (with respect to either orientation).
\item The \emph{trivial circles} have zero $\sltwo$ weight space  with respect to $*$ (with respect to either orientation).
\ei
The ($\sltwo$ weight space) \emph{$k$-grading} of $ v$ is given by 
$$ 
gr_k(v) = \#\{ p_i :\  p_i = +,\  C_i \text{ nontrivial}\} 
	- \# \{ p_i:\  p_i = -, \ C_i \text{ nontrivial}\}.
	$$

\ei

We think of the $i$-th tensor factor in $\mathbb{V}^{\otimes |R_\alpha|}$ as corresponding to a labeling of the $i$-th component $C_i$ with a $-$ or $+$ symbol. As such, for a trivial circle $C_i$ we will henceforth write the corresponding factor as $\mathbb{W} = \langle w_-, w_+\rangle$, reserving the notation $\mathbb{V} = \langle v_-, v_+ \rangle$ for nontrivial circles.
For example, if $R_\alpha$ consists of one nontrivial circle and two trivial circles, labeled in that order, one of the generators of $V_\alpha$ is $v_- \otimes w_+ \otimes w_-$. 

\brmk
In \cite{GLW-2015}, Grigsby, Licata, and Wehrli show that annular Khovanov homology over $\C$ enjoys a Lie algebra $\sltwo(\C)$-action, which is why we refer to the $k$-grading as the ``$\sltwo$ weight space grading.''
\ermk

\subsubsection*{Edges} \label{edges}

Let $\alpha, \alpha' \in \{ 0, 1 \}^n$.  If $\alpha'$ is a bit increment from $\alpha$, we say that the resolution $R_\alpha'$ is a \emph{successor to} $R_\alpha$. In our cube of resolutions, we draw a directed edge $R_\alpha \to R_\alpha'$ corresponding to a linear map 
$\partial_{\alpha, \alpha'}: V_\alpha \to V_{\alpha'}$.
This map depends on a few characteristics of the resolutions $R_\alpha$ and $R_{\alpha'}$. 

\bi

\item Since $R_\alpha$ and $R_{\alpha'}$ differ at the site of exactly one crossing, all but one or two components remain unchanged between $R_\alpha$ and $R_{\alpha'}$. 
The linear map $\partial_{\alpha, \alpha'}$ may be a \emph{merge map} or a \emph{split map}.
If the site of the crossing touches only one component $C_i$ of $R_\alpha$, then in $R_{\alpha'}$, $C_i$ has split into two components. 
If the site of the crossing touches two components $C_i$ and $C_j$ of $R_\alpha$, then in $R_{\alpha'}$, $C_i$ and $C_j$ merge into one component.

\item The participating components may be trivial or nontrivial circles. The linear map $\partial_{\alpha, \alpha'}$ is constructed so that $gr_j$ and $gr_k$ are preserved, and $gr_i$ is increased by 1.
Roberts shows in \cite{Roberts-2013} that there are exactly six possible types of maps corresponding to an edge in the cube resolutions.
We can visually describe these maps by drawing dots representing the distinguished generators of $V_\alpha$ and $V_\alpha'$, aligning them in rows of constant $gr_j$, and drawing arrows between the dots to represent components of the differential. Figure \ref{akh-maps} shows the six maps defining all possible $\partial_{\alpha, \alpha'}$ using this visual respresentation.

\ei

\begin{figure}
\centering

\begin{minipage}{0.47\textwidth}
	\centering
	\begin{tikzpicture}
	\filldraw 
	(0,0) circle (2pt) node[align=center, below = 0.1 cm]{$v_-$} 
	(0,2) circle (2pt) node[align = center, below = 0.1 cm]{$v_+$}
	(4,4) circle (2pt) node[align=center, below = 0.1 cm]{$v_+\otimes w_+$}
	(3,2) circle (2pt) node[align=center, below = 0.1 cm]{$v_+ \otimes w_-$} 
	(5,2) circle (2pt) node[align=center, below = 0.1 cm]{$v_- \otimes w_+$} 
	(4,0) circle (2pt) node[align=center, below = 0.1 cm]{$v_- \otimes w_-$};
	\path[->] (0.2,0) edge (3.8,0);
	\path[->] (0.2,2) edge (2.8, 2);
	\end{tikzpicture}
	\caption*{Type A: $v \to vw$}
	\label{v-vw}
\end{minipage} 
\begin{minipage}{0.47\textwidth}
	\centering
\begin{tikzpicture}
\filldraw 
(8,2) circle (2pt) node[align=center, below = 0.1 cm]{$v_-$} 
(8,4) circle (2pt) node[align = center, below = 0.1 cm]{$v_+$}
(4,4) circle (2pt) node[align=center, below = 0.1 cm]{$v_+\otimes w_+$}
(3,2) circle (2pt) node[align=center, below = 0.1 cm]{$v_+ \otimes w_-$} 
(5,2) circle (2pt) node[align=center, below = 0.1 cm]{$v_- \otimes w_+$} 
(4,0) circle (2pt) node[align=center, below = 0.1 cm]{$v_- \otimes w_-$};
\path[->] (5.2,2) edge (7.8, 2);
\path[->] (4.2, 4) edge (7.8, 4);
\end{tikzpicture} 
	\caption*{Type D: $vw \to v$}
	\label{vw-v}
\end{minipage}

\begin{minipage}{0.47\textwidth}
	\centering
	\vspace{.5in}
	\begin{tikzpicture}
	\filldraw 
	(0,0) circle (2pt) node[align=center, below = 0.1 cm]{$w_-$} 
	(0,2) circle (2pt) node[align = center, below = 0.1 cm]{$w_+$}
	(4,4) circle (2pt) node[align=center, below = 0.1 cm]{$v_+\otimes v_+$}
	(3,2) circle (2pt) node[align=center, below = 0.1 cm]{$v_+ \otimes v_-$} 
	(5,2) circle (2pt) node[align=center, below = 0.1 cm]{$v_- \otimes v_+$} 
	(4,0) circle (2pt) node[align=center, below = 0.1 cm]{$v_- \otimes v_-$};
	\path[->] (0.2,2) edge[bend left] (2.8, 2);
	\path[->] (0.2,2) edge[bend left] (4.8,2);
	\end{tikzpicture}
	\caption*{Type B: $w \to vv$}
\end{minipage}
\begin{minipage}{0.47\textwidth}
	\centering
	\vspace{.5in}
\begin{tikzpicture}
\filldraw 
(8,2) circle (2pt) node[align=center, below = 0.1 cm]{$w_-$} 
(8,4) circle (2pt) node[align = center, below = 0.1 cm]{$w_+$}
(4,4) circle (2pt) node[align=center, below = 0.1 cm]{$v_+\otimes v_+$}
(3,2) circle (2pt) node[align=center, below = 0.1 cm]{$v_+ \otimes v_-$} 
(5,2) circle (2pt) node[align=center, below = 0.1 cm]{$v_- \otimes v_+$} 
(4,0) circle (2pt) node[align=center, below = 0.1 cm]{$v_- \otimes v_-$};
\path[->] (3.2,2) edge[bend left] (7.8, 2);
\path[->] (5.2,2) edge[bend left] (7.8, 2);
\end{tikzpicture}
	\caption*{Type E: $vv \to w$}
	\label{vv-w}
\end{minipage} 

\begin{minipage}{0.47\textwidth}
	\centering
\vspace{.5in}
\begin{tikzpicture}
\filldraw 
(0,0) circle (2pt) node[align=center, below = 0.1 cm]{$w_-$} 
(0,2) circle (2pt) node[align = center, below = 0.1 cm]{$w_+$}
(4,4) circle (2pt) node[align=center, below = 0.1 cm]{$w_+\otimes w_+$}
(3,2) circle (2pt) node[align=center, below = 0.1 cm]{$w_+ \otimes w_-$} 
(5,2) circle (2pt) node[align=center, below = 0.1 cm]{$w_- \otimes w_+$} 
(4,0) circle (2pt) node[align=center, below = 0.1 cm]{$w_- \otimes w_-$};
\path[->] (0.2,0) edge (3.8,0);
\path[->] (0.2,2) edge[bend left] (2.8, 2);
\path[->] (0.2,2) edge[bend left] (4.8,2);
\end{tikzpicture}
	\label{w-ww}
	\caption*{Type C: $w \to ww$}
\end{minipage}
\begin{minipage}{0.47\textwidth}
	\centering
	\vspace{.5in}
\begin{tikzpicture}
\filldraw 
(8,2) circle (2pt) node[align=center, below = 0.1 cm]{$w_-$} 
(8,4) circle (2pt) node[align = center, below = 0.1 cm]{$w_+$}
(4,4) circle (2pt) node[align=center, below = 0.1 cm]{$w_+\otimes w_+$}
(3,2) circle (2pt) node[align=center, below = 0.1 cm]{$w_+ \otimes w_-$} 
(5,2) circle (2pt) node[align=center, below = 0.1 cm]{$w_- \otimes w_+$} 
(4,0) circle (2pt) node[align=center, below = 0.1 cm]{$w_- \otimes w_-$};
\path[->] (4.2,4) edge (7.8,4);
\path[->] (3.2,2) edge[bend left] (7.8, 2);
\path[->] (5.2,2) edge[bend left] (7.8,2);
\end{tikzpicture}
	\caption*{Type F: $ww \to w$}
	\label{v-vw}
\end{minipage}

\caption{The $AKh$ differentials. Within each diagram, the rows have constant $j$-grading.}
\label{akh-maps}
\end{figure}

One can check that $\partial_{AKh} = \sum \partial_{\alpha, \alpha'}$ is a differential by verifying that each face of the cube of chains commutes.
The triply graded homology of this complex $CKh(\CD(L))$ is independent of the choice of diagram $\CD(L)$ \cite{APS-2004}.

\section{The annular Khovanov-Tate bicomplex} \label{annular-khovanov-tate-bicomplex}

We are now prepared to construct the bicomplex which defines the spectral sequence in Theorem \ref{main-thm}. The motivation for this construction is discussed in \S \ref{equiv-hom-tate-ss}.

\subsection{An involution on $CKh$}  \label{involution-on-ckh}

The first task is to define an involution on $CKh(\CD(\tilde L))$ induced by the topological involution $\tau$.

\bdefn \label{defn-taus}
We will need to discuss many related involutions. By abuse of notation, we utilize only three symbols; the distinctions within each group will be clear from context. After presenting the notation, we'll discuss the involutions and their relationships with each other in more detail.

In the following contexts, we use the symbol $\tau$:
\bi
\item the bitstring involution defined in Definition \ref{bitstring-notation}
\item the topological involution on $S^3$ with fixed point set $B$
\item the restriction of $\tau$ to the 2-periodic link $\tilde L \subset S^3$
\item the involution induced on $A \times I$ and its restriction to $\tilde L \hookrightarrow A \times I$.
\ei

We use the symbol $\tau_A$ to denote the following:
\bi
\item the involution induced on $A$ from $\tau$ on $A \times I$
\item the restriction of $\tau_A$ to $\CD(\tilde L) \subset A$
\item the involution induced on an equivariant resolution $R_{\alpha\alpha} \subset A$.
\ei

Finally, $\tau^\#$ is used to denote the following:
\bi
\item the involution induced on the chains of $CKh(\CD(\tilde L))$
\item the involutions induced by $\tau^\#$ on the pages of the spectral sequences.
\ei

\edefn

We continue to use the setting described in \S \ref{2-periodic-links-diagrams}. 
Let $n$ be the number of crossings in the quotient link diagram $\CD(L)$. 
Label the $2n$ crossings in $\CD(\tilde L)$ by first labeling crossings of $\CD(\tilde L)$ on the right side of $A$ ($x > 0$) , then labeling the crossings on the left side ($x < 0$) so that $\tau_A$ takes the $i$-th crossing on the right to the $(i+n)$-th crossing on the left.

With this assignment, each $2n$-bit string $\alpha = \alpha_1\alpha_2$ is the concatenation of two $n$-bit strings $\alpha_1$ and $\alpha_2$, where $\alpha_1$ represents a sequence of resolution choices for the crossings on the right and $\alpha_2$ represents a sequence of resolution choices for the crossings on the left. 
Thus $\tau_A$ takes $R_\alpha = R_{\alpha_1\alpha_2}$ to $R_{\tau(\alpha)} = R_{\alpha_2\alpha_1}$.

Thinking of the distinguished generators of the cube of chains $\scube(\CD(\tilde L))$ as marked resolutions, $\tau_A$ induces an involution $\tau^\#$ on the chains of $CKh(\CD(\tilde L))$. 

\bdefn \label{defn-equivariant-resolutions}
A resolution $R_\alpha$ is \emph{equivariant} if $\tau_A(R_\alpha) = R_\alpha$.
A generator $x$ of the cube of chains $\scube(\CD(\tilde L))$ is \emph{equivariant} if $\tau^\#(x) = x$.
If a resolution or generator is not equivariant, it is \emph{nonequivariant}.
\edefn

In fact, $\tau^\#$ is an involution on the complex $(CKh(\CD(\tilde L)), \partial_{AKh})$:

\blem \label{akh-and-tau-commute}

$\partial_{AKh}$ and $\tau^\#$ commute.

\begin{proof}

It suffices to show that every edge belongs to a well-defined pair of edges $\{\partial_{\alpha, \alpha'}, \partial_{\tau(\alpha), \tau(\alpha')}\}$ of $\text{\sc Cube}(\CD(\tilde L))$, where $\alpha'$ is a bit increment from $\alpha$, so that the following diagram commutes.
\begin{center}
\begin{tikzcd}
V_{\alpha'} \arrow{r}{\tau^\#} & \tau^\#(V_{\alpha'}) \\
V_{\alpha} \arrow{r}{\tau^\#}  \arrow{u}{\partial_{\alpha, \alpha'}} & \tau^\#(V_\alpha) 
	\arrow[swap]{u}{\partial_{\tau(\alpha), \tau(\alpha')}}
\end{tikzcd}
\end{center}
Since $\partial_{AKh}$ is the sum of all such pairs of edges, $\partial_{AKh}$ commutes with $\tau^\#$.

The edge $\partial_{\alpha, \alpha'}$ may correspond to a change of resolution on the right or left side of $\CD(\tilde L)$. Without loss of generality, consider the case where the resolution change is on the right.
\begin{center}
\begin{tikzcd}
V_{\alpha'_1\alpha_2} \arrow{r}{\tau^\#} & V_{\alpha_2\alpha_1'} \\
V_{\alpha_1\alpha_2} \arrow{r}{\tau^\#}  \arrow{u}{\partial_{\alpha, \alpha'}} & 													V_{\alpha_2\alpha_1}	\arrow[swap]{u}{\partial_{\tau(\alpha),\tau(\alpha')}}
\end{tikzcd}
\end{center}
The components of $R_{\alpha_1\alpha_2}$ and $R_{\alpha_2\alpha_1}$ are identified in pairs by the diffeomorphism $\tau_A$; the same holds for the pair $R_{\alpha_1'\alpha_2}$ and $R_{\alpha_2\alpha_1'}$.
In this way, a merge (resp.\ split) map $V_{\alpha_1\alpha_2} \to V_{\alpha_1'\alpha_2}$ corresponds to the isomorphic merge (resp.\ split) map $V_{\alpha_2 \alpha_1} \to V_{\alpha_2\alpha_1'}$. 

Notice that $ \partial_{\alpha, \alpha'} = \partial_{\tau(\alpha)}$ if and only if they have the same source and target. This is impossible as one of $|\alpha|$ and $|\alpha'|$ is odd, and hence corresponds to a nonequivariant resolution.
\end{proof}

\elem

\subsection{Construction of the $\atate$ bicomplex} \label{construction-of-the-atate-bicomplex}

This allows us to define a variant of the Tate bicomplex (see \S \ref{equiv-hom-tate-ss}), constructed as follows:

\bi
\item Each column is a copy of the complex $CKh(\CD(\tilde L))$, where the vertical filtration is induced by $gr_i$.
\item The underlying vector space is the $\F_2[\theta]$-module $CKh(\CD(\tilde L)) \otimes_{\F_2} \F_2[\theta]$, where $\theta$ acts on the columns by shifting right by one column.
\item The horizontal differential is $1 + \tau^\#$.
\ei

Observe that this complex is quadruply graded by $gr_i, gr_j, gr_k,$ and a column grading $gr_t$, where $t$ is the exponent of $\theta$.

\bdefn \label{defn-annular-khovanov-tate-complex}
We call this bicomplex the \emph{annular Khovanov-Tate complex} for the 2-periodic link $\tilde L$ with involution $\tau$, abbreviated $\atate(\tilde L)$. 
We denote the total complex of this bicomplex by $Tot(\atate(\tilde L))$. 
\edefn

\brmk
This notation requires some justification. We'll show that 
$$H^*(Tot(\atate(\tilde L))) \cong AKh(L)^{2j-k,k} \otimes \F_2[\theta, \theta\inv]$$
 where the latter is independent of the choice of diagram by \cite{APS-2004}.
\ermk

\section{Proof of Theorem \ref{main-thm}} \label{proof-of-theorem}

We are now ready to prove the main result. Recall that $AKh^{j,k}(L)$ denotes the part of $AKh(L)$ at quantum grading $j$ and $\sltwo$ weight space  grading $k$.

\begin{reptheorem}{main-thm}
Let $\tilde L$ be a 2-periodic link with quotient link $L$. For each pair of integers $(j,k)$, there is a spectral sequence with
$$
E^1 \cong AKh^{2j-k,k}(\tilde L) \otimes_{\FF} \FF[\theta, \theta\inv] 
\abuts
E^\infty \cong AKh^{j,k}(L) \otimes_{\FF} \FF[\theta, \theta\inv].
$$
\end{reptheorem}

Before diving into the details, we first give a sketch of the main ideas behind the proof. The reader may also wish to refer to the example in \S \ref{eg-hopf-link}. 

Since the complex  $CKh(\CD(\tilde L))$ is finite-dimensional, our bicomplex $\atate(\CD(\tilde L))$ has finite-dimensional columns, so by Proposition \ref{hvE-and-vhE}, the $\hvE$ and $\vhE$ spectral sequences both converge to $H_*(Tot(\atate(\CD(\tilde L))))$. 
We will compute the $\hvE$ spectral sequence to find that $\hvE^\infty \cong AKh(L) \otimes \F_2[\theta, \theta\inv]$. Then we'll show that $\vhE^1 = AKh(\tilde L)\otimes \F_2[\theta, \theta\inv]$. By Proposition \ref{hvE-and-vhE}, $\vhE^\infty \cong \hvE^\infty$, with the isomorphism respecting the $(i+t)$-grading. Observing that each diagonal level set in $\vhE^\infty$ can be identified with $\oplus_i AKh_i(L)$ (where $AKh_i(L)$ refers to the piece of $AKh(L)$ at $gr_i = i$), we see that $\vhE$ is a spectral sequence from the $(i+t)$-filtered $AKh(\tilde L) \otimes \F_2[\theta, \theta\inv]$ to $AKh(L) \otimes \F_2[\theta, \theta\inv]$. The total differential preserves $j$- and $k$- gradings, so $\vhE$ splits along $j$- and $k$- gradings.

\subsection{The first pages} \label{the-first-pages}

The first pages of either spectral sequence can be understood readily and have topological interpretations.

\blem \label{vhE1}
$\vhE^1 \cong AKh(\tilde L) \otimes_{\FF} \FF[\theta, \theta\inv]$.

\bpf

In each column, the cancellation of all the vertical arrows computes the homology $AKh(\CD(\tilde L))$. 
The fact that cancellation actually occurs in the greater context of the bicomplex is irrelevant because any induced maps are no longer vertical.
\epf
\elem

\blem \label{hvE1-gen-by-equiv-gens}
$\hvE^1$ is generated by the equivariant generators of $\atate(\CD(\tilde L))$.

\begin{figure} 
\centering 
\def\svgwidth{400pt} 
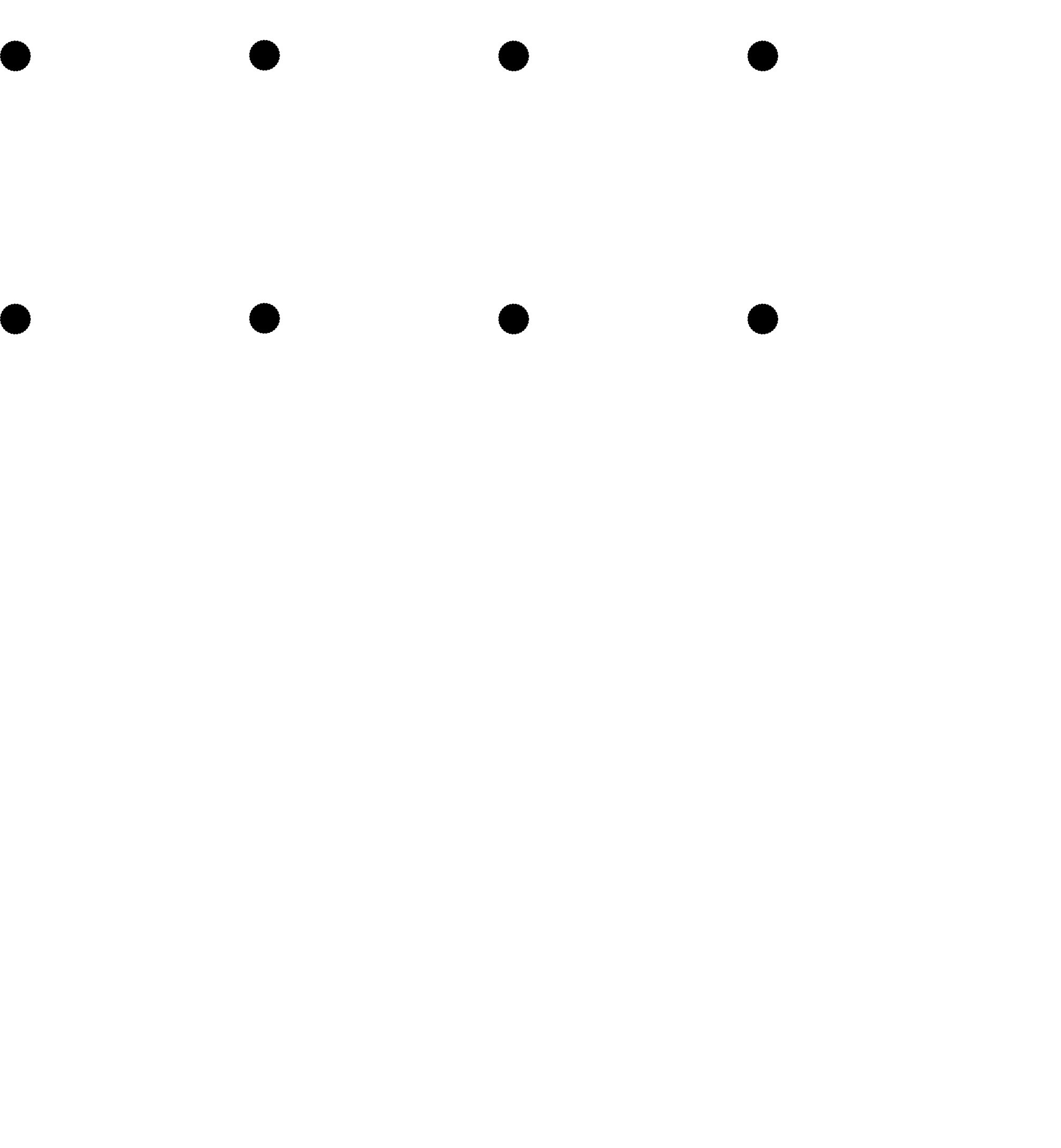
\caption{Cancellation of the $\tau^\#$ arrows eliminates all the nonequivariant generators.}
\label{cancel-all-taus} 
\end{figure}

\bpf
Since $1+\tau^\#$ vanishes on all equivariant generators, it suffices to show that all the nonequivariant generators vanish upon cancellation of all the horizontal arrows in the bicomplex.
In fact, cancellation of all the $\tau^\#$ arrows suffices, as shown in Figure \ref{cancel-all-taus}.

\epf
\elem

\bdefn  \label{defn-staircase-of-arrows}

\begin{figure} 
\centering 
\def\svgwidth{150pt} 
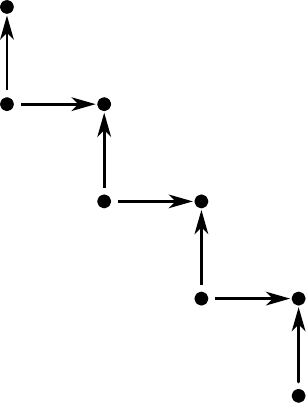
\caption{A staircase of arrows. In this example, the staircase has \emph{height} 4, so it induces an arrow of \emph{length} 4.}
\label{staircase} 
\end{figure}

A \emph{staircase of arrows} (Figure \ref{staircase}) is a finite collection of arrows in a bicomplex (generated by a distinguished filtered graded basis) which form a path with the following properties.
\bi
\item The path alternates between vertical and horizontal arrows.
\item Both arrows at a corner either both point away from or both point toward the corner.
\item If a vertical arrow and a horizontal arrow meet at a corner of the path, they either have the same source or the same target.
\ei

The \emph{height} of a staircase of arrows with respect to a given filtration is the change in filtration degree between the generators at the ends of the path.

With respect to the row-wise filtration, if both ends of the staircase begin with a vertical arrow, then cancellation of the horizontal arrows in the staircase (in any order) induces a single arrow from one end of the path to the other. The \emph{length} of this induced arrow is the height of the staircase.
\edefn

\blem \label{dr-odd-zero}
If $r$ is odd, $^{hv}d^r = 0$.
\bpf
By Lemma \ref{hvE1-gen-by-equiv-gens}, only equivariant generators survive to $\hvE^1$. Since these belong to vector spaces corresponding to equivariant resolutions of $\CD(\tilde L)$, their Hamming weights must be even. The $i$-grading is a shift of Hamming weight, so every other row vanishes on $\hvE^1$. Components of $^{hv}d_r$ for odd $r$ are represented by arrows of length $r$, which must map to or from a vanishing row, so $^{hv}d_r = 0$. 
\epf
\elem

Therefore the next interesting differential in the spectral sequence is $^{hv}d^2$.

\blem \label{d2-from-cancellation}
Every component of $^{hv} d^2$ is induced by the cancellation of a $\tau^\#$ arrow at a row corresponding to an odd Hamming weight.

\bpf
Since we begin with only (row-wise) filtration degree 0 and 1 maps, the only way to induce a length 2 differential is by canceling the horizontal arrow on a staircase of height 2. 
The source and target of a length 2 differential is an equivariant generator, so the cancelled horizontal arrow must lie in a row corresponding to an odd Hamming weight.
\epf
\elem

\brmk \label{cancellation-and-lift}
Later on it will become useful to view the computation of the induced differentials in the spectral sequence as a sequence of maps (by $\partial_{AKh}$) and lifts (by $1 + \tau^\#$), traveling up the staircase in Figure \ref{staircase}. The maps constructed this way are the same as the maps induced by cancellation of the $\tau^\#$ arrows in the staircase. This is a more intuitive reason why cancellation of only the $\tau^\#$ arrows is sufficient for eliminating the nonequivariant generators.
\ermk

\subsection{Relationship between $(\hvE^2, ^{hv}d^2)$ and the quotient link}

The crux of the theorem comes from the observation that the generators and differentials on page 2 of the $\hvE$ spectral sequence correspond exactly with distinguished generators in the cube of chains for the diagram of the quotient link. The following two propositions verify the details of this correspondence.

\bprop \label{E2-correspondence}
There is a one-to-one correspondence

\begin{align*}
\left \{ \text{generators of } \hvE^2 \right \}  \longleftrightarrow &
\left \{ \text{generators of } CKh(L) \otimes_{\FF} \FF[\theta, \theta\inv] \right \}\\
\tilde v  \longleftrightarrow &  v
\end{align*}

induced by $\tau$ such that 
\bi
\item $gr_i(\tilde v) = 2 gr_i(v)$
\item $gr_j(\tilde v) = 2 gr_j(v) - gr_k(v)$
\item $gr_k(\tilde v) = gr_k(v)$
\item $gr_t(\tilde v) = gr_t(v)$.
\ei

\bpf
By Lemma \ref{hvE1-gen-by-equiv-gens}, $\hvE^1$ is generated by the equivariant generators of $\atate(\CD(\tilde L))$, and Lemma \ref{dr-odd-zero} shows $^{hv}d^1 =0$, so $\hvE^1 =\hvE^2$.

\subsubsection*{Bijection of generators}
First of all, the correspondence is given by the identity on the $\FF[\theta, \theta\inv]$ factor. 

Let $\tilde v$ denote a generator of $\hvE^2$, coming from the equivariant resolution $R_{\alpha\alpha}$.  Thinking of $\tilde v$ as a $\tau_A$-equivariant assignment of plusses and minuses to the components of $R_{\alpha\alpha}$, let $v$ be the quotient of $\tilde v$ by $\tau_A$. 

In the other direction, for $v$ a generator of $CKh(L) \otimes_{\FF} \FF[\theta, \theta\inv]$ coming from a resolution $R_{\alpha}$, view $v$ as an assignment of plusses and minuses to the components of $R_\alpha$. Pick any path connecting the two components of $\partial A$ and take a double cover of $R_\alpha$. The lift of $v$ is $\tilde v$.

\subsubsection*{Grading relationships}
The $t$-grading relationship follows by definition. For the other three relationships, suppose $v$ consists of 
\bi
\item $a$ nontrivial circles labeled ``$+$''
\item $b$ nontrivial circles labeled ``$-$''
\item $c$ trivial circles labeled ``$+$''
\item $d$ trivial circles labeled ``$-$''
\ei
In this case, write $v = v_+^a v_-^b w_+^c w_-^d$ to denote the distinguished generator, suppressing the tensor signs in favor of compact notation.
Then $\tilde v = v_+^a v_-^b w_+^{2c} w_-^{2d}$. 

Let $n_+$ and $n_-$ denote the number of positive and negative crossings in $\CD(L)$, respectively. 
For $v$, 
\bi
\item $gr_i(v) = |\alpha| - n_-$
\item $gr_j(v) = |\alpha| + (a-b+c-d) +n_+ - 2 n_-$
\item $gr_k(v) = a - b$.
\ei

Since the diagram $\CD(\tilde L)$ obtained by taking the double cover of $\CD(L)$ has twice as many positive and negative crossings as $\CD(L)$, respectively, we compute that
\bi
\item $gr_i(\tilde v) = 2|\alpha| - 2n_- = 2 gr_i(v)$
\item $gr_j(\tilde v) = 2|\alpha| + (a - b + 2c - 2d) +2n_+ - 4n_- 
= 2 gr_j(v) - (a-b) = 2 gr_j(v) - gr_k(v)$
\item $gr_k(\tilde v) = a-b= gr_k(v)$.
\ei
\epf
\eprop

\bprop \label{d2-correspondence}
Under the correspondence above, each slope -2 line in the bicomplex $(\hvE^2, ^{hv} d^2)$ is isomorphic to the complex $(CKh(L), \partial_{AKh})$.

\bpf
Let $\alpha'$ be a bit increment from $\alpha$. 
Lemma \ref{d2-from-cancellation} implies that any component of $^{hv}d^2$ lies along a slope -2 line. 

\begin{figure} 
\centering 
\def\svgwidth{150pt} 
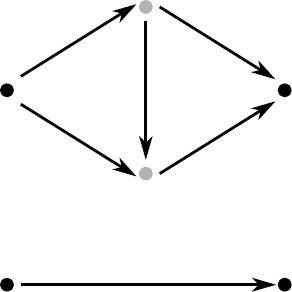
\caption{Cancellation of $\tau^\#$ arrows induce $\dtwo$ arrows. For each diamond of the form above, we obtain a length 2 differential.}
\label{diamond} 
\end{figure}

Upstairs, consider the components of $^{hv}d^2$ from $V_{\alpha\alpha} \otimes \theta^t$ to $V_{\alpha'\alpha'} \otimes \theta^{t-1}$. (See Figure \ref{diamond}.)
The two resolutions $R_{\alpha\alpha'}$ and $R_{\alpha'\alpha}$ connecting $R_{\alpha\alpha}$ and $R_{\alpha'\alpha'}$ correspond  to changing the resolutions of some pair of crossings $c_i$ and $c_{i+n}$ from their 0-resolution in $R_{\alpha\alpha}$ to their 1-resolution in $R_{\alpha'\alpha'}$.

Downstairs, consider the edge of $\scube(\CD(L))$ going from $V_\alpha$ to $V_{\alpha'}$. This edge falls into one of the six types of $AKh$ differentials, as discussed in Section \ref{annular-khovanov-homology}. 

We compare the $^{hv}d^2$ components induced by the cancellation of the $\tau^\#$ arrows going between $V_{\alpha\alpha'}$ and $V_{\alpha'\alpha}$ with the components of $\partial_{AKh}$, in the six different cases. Figures \ref{v-vww} through \ref{ww-wwww} show the computation of $\dtwo$ in all six cases. 
Compare these figures with Figure \ref{akh-maps}. Note that in contrast with the notation used in Figure \ref{akh-maps}, the tensor signs are now suppressed for the sake of compactness.
By inspection, under the correspondence from Proposition \ref{E2-correspondence}, the induced length 2 arrows upstairs align exactly with the differential downstairs.


\bfig
\begin{minipage}{1.2\textwidth}
\centering
\includegraphics[height=2in]{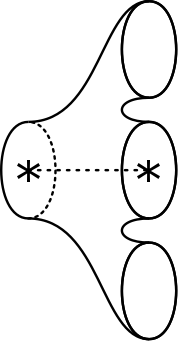}
	\begin{tikzpicture}
	\filldraw 
	(0,0) circle (3pt) node[align=center, below=0.1 cm]{$v_-$}
	(0,2) circle (3pt) node[align=center, below=0.1 cm]{$v_+$}
	(2,2) circle (1pt) node[align=center, below=0.1 cm]{$v_+ w_-$}
	(3,0) circle (1pt) node[align=center, below=0.1 cm]{$v_-  w_-$}
	(3,4) circle (1pt) node[align=center, below=0.1 cm]{$v_+  w_+$}
	(4,2) circle (1pt) node[align=center, below=0.1 cm]{$v_- w_+$}
	(6,2) circle (3pt) node[align=center, below=0.1 cm]{$v_+ w_-  w_-$}
	(6,4) circle (1pt) node[align=center, below=0.1 cm]{$v_+  w_+  w_-$}
	(8,0) circle (3pt) node[align=center, below=0.1 cm]{$v_-  w_-  w_-$}
	(8,2) circle (1pt) node[align=center, below=0.1 cm]{$v_-  w_+  w_-$}
	(8,4) circle (1pt) node[align=center, below=0.1 cm]{$v_+  w_-  w_+$}
	(8,6) circle (3pt) node[align=center, below=0.1 cm]{$v_+  w_+  w_+$}
	(10,2) circle (1pt) node[align=center, below=0.1 cm]{$v_- w_-  w_+$}
	(10,4) circle (3pt) node[align=center, below=0.1 cm]{$v_-  w_+  w_+$};
	\path[->] (0.2,2) edge[bend left] (1.8,2);
	\path[->] (0.2,0) edge (2.8,0);
	\path[->] (2.2,2) edge[bend left] (5.8,2);
	\path[->] (3.2,0) edge (7.8,0);
	\end{tikzpicture}
	\caption{Type $\mathscr{A}$: $\tilde v \to \tilde v \tilde w$}
	\label{v-vww}
\end{minipage} 

\vspace{0.5 in}
\begin{minipage}{1.2\textwidth}
\centering
\includegraphics[height=2in]{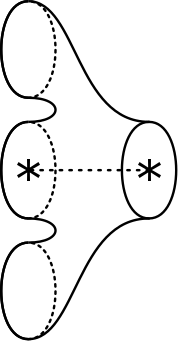}
	\begin{tikzpicture}
	\filldraw 
	(16,4) circle (3pt) node[align=center, below=0.1 cm]{$v_-$}
	(16,6) circle (3pt) node[align=center, below=0.1 cm]{$v_+$}
	(12,4) circle (1pt) node[align=center, below=0.1 cm]{$v_+ w_-$}
	(13,2) circle (1pt) node[align=center, below=0.1 cm]{$v_-  w_-$}
	(13,6) circle (1pt) node[align=center, below=0.1 cm]{$v_+  w_+$}
	(14,4) circle (1pt) node[align=center, below=0.1 cm]{$v_- w_+$}
	(6,2) circle (3pt) node[align=center, below=0.1 cm]{$v_+ w_-  w_-$}
	(6,4) circle (1pt) node[align=center, below=0.1 cm]{$v_+  w_+  w_-$}
	(8,0) circle (3pt) node[align=center, below=0.1 cm]{$v_-  w_-  w_-$}
	(8,2) circle (1pt) node[align=center, below=0.1 cm]{$v_-  w_+  w_-$}
	(8,4) circle (1pt) node[align=center, below=0.1 cm]{$v_+  w_-  w_+$}
	(8,6) circle (3pt) node[align=center, below=0.1 cm]{$v_+  w_+  w_+$}
	(10,2) circle (1pt) node[align=center, below=0.1 cm]{$v_- w_-  w_+$}
	(10,4) circle (3pt) node[align=center, below=0.1 cm]{$v_-  w_+  w_+$};
	\path[->] (8.2,6) edge (12.8,6);
	\path[->] (13.2,6) edge (15.8,6);
	\path[->] (10.2,4) edge[bend left] (13.8,4);
	\path[->] (14.2,4) edge[bend left] (15.8,4);
	\path[->] (6.2,2) edge[bend left] (12.8,2);
	\end{tikzpicture}
	\caption{Type $\mathscr{D}$: $ \tilde v \tilde w \to \tilde v$}
	\label{vww-v}
\end{minipage} 
\efig

\bfig
\begin{minipage}{1.2\textwidth}
\centering
\includegraphics[height=2in]{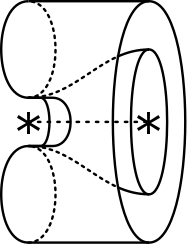}
	\begin{tikzpicture}
	\filldraw 
	(0,0) circle (3pt) node[align=center, below=0.1 cm]{$w_-w_-$}
	(-1,2) circle (1pt) node[align=center, below=0.1 cm]{$w_+w_-$}
	(1,2) circle (1pt) node[align=center, below=0.1 cm]{$w_-w_+$}
	(0,4) circle (3pt) node[align=center, below=0.1 cm]{$w_+w_+$}
	(3,2) circle (1pt) node[align=center, below=0.1 cm]{$w_-$}
	(3,4) circle (1pt) node[align=center, below=0.1 cm]{$w_+$}
	(5,4) circle (3pt) node[align=center, below=0.1 cm]{$v_+v_-$}
	(6,2) circle (3pt) node[align=center, below=0.1 cm]{$v_-v_-$}
	(6,6) circle (3pt) node[align=center, below=0.1 cm]{$v_+v_+$}
	(7,4) circle (3pt) node[align=center, below=0.1 cm]{$v_-v_+$};
	\path[->] (0.2,4) edge (2.8,4);
	\path[->] (3.2,4) edge[bend left] (4.8,4);
	\path[->] (3.2,4) edge[bend left] (6.8,4);
	\end{tikzpicture}
	\caption{Type $\mathscr{B}$: $\tilde w \to \tilde v \tilde v$}
	\label{ww-vv}
\end{minipage}

\begin{minipage}{1.2\textwidth}
\centering
\hspace{-1in}\includegraphics[height=2in]{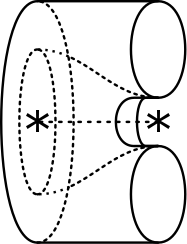}
	\begin{tikzpicture}
	\filldraw 
	(7,4) circle (3pt) node[align=center, below=0.1 cm]{$w_-w_-$}
	(6,6) circle (1pt) node[align=center, below=0.1 cm]{$w_+w_-$}
	(8,6) circle (1pt) node[align=center, below=0.1 cm]{$w_-w_+$}
	(7,8) circle (3pt) node[align=center, below=0.1 cm]{$w_+w_+$}
	(4,4) circle (1pt) node[align=center, below=0.1 cm]{$w_-$}
	(4,6) circle (1pt) node[align=center, below=0.1 cm]{$w_+$}
	(0,4) circle (3pt) node[align=center, below=0.1 cm]{$v_+v_-$}
	(1,2) circle (3pt) node[align=center, below=0.1 cm]{$v_-v_-$}
	(1,6) circle (3pt) node[align=center, below=0.1 cm]{$v_+v_+$}
	(2,4) circle (3pt) node[align=center, below=0.1 cm]{$v_-v_+$};
	\path[->] (0.2,4) edge[bend left] (3.8,4);
	\path[->] (4.2,4) edge (6.8,4);
	\path[->] (2.2,4) edge[bend left] (3.8,4);
	\end{tikzpicture}
	\caption{Type $\mathscr{E}$: $\tilde v \tilde v \to \tilde w$}
	\label{vv-ww}
\end{minipage} 
\efig

\bfig
\begin{minipage}{1.2\textwidth}
\centering

\hspace{-1in} \includegraphics[height=2in]{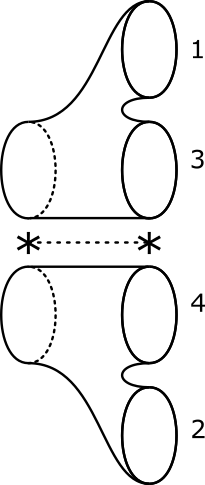}

\vspace{0.5in}

\hspace{-1in}
	\begin{tikzpicture}
	\filldraw 
	(0,2) circle (1pt) node[align=center, below=0.1 cm]{$w_+ w_-$}
	(1,0) circle (3pt) node[align=center, below=0.1 cm]{$w_- w_-$}
	(2,2) circle (1pt) node[align=center, below=0.1 cm]{$w_-w_+$}
	(1,4) circle (3pt) node[align=center, below=0.1 cm]{$w_+w_+$}
	(4,2) circle (1pt) node[align=center, below=0.1 cm]{$w_+w_-w_-$}
	(4,4) circle (1pt) node[align=center, below=0.1 cm]{$w_+w_+w_-$}
	(6,0) circle (1pt) node[align=center, below=0.1 cm]{$w_-w_-w_-$}
	(6,2) circle (1pt) node[align=center, below=0.1 cm]{$w_-w_+w_-$}
	(6,4) circle (1pt) node[align=center, below=0.1 cm]{$w_+w_-w_+$}
	(6,6) circle (1pt) node[align=center, below=0.1 cm]{$w_+w_+w_+$}
	(8,2) circle (1pt) node[align=center, below=0.1 cm]{$w_-w_-w_+$}
	(8,4) circle (1pt) node[align=center, below=0.1 cm]{$w_-w_+w_+$}
	(10,4)  circle (3pt) node[align=center, below=0.1 cm]{$w_+w_+ w_-w_-$}
	(12.5,0) circle (3pt) node[align=center, below=0.1 cm]{$w_- w_-w_- w_-$}
	(12.5,8) circle (3pt) node[align=center, below=0.1 cm]{$w_+w_+w_+w_+$}
	(15,4) circle (3pt) node[align=center, below=0.1 cm]{$w_- w_-w_+w_+$}
	(11,2) circle (1pt) node{}
	(11,4) circle (1pt) node{}
	(11,6) circle (1pt) node{}
	(12,2) circle (1pt) node{}
	(12,4) circle (1pt) node{}
	(12,6) circle (1pt) node{}
	(13,2) circle (1pt) node{}
	(13,4) circle (1pt) node{}
	(13,6) circle (1pt) node{}
	(14,2) circle (1pt) node{}
	(14,4) circle (1pt) node{}
	(14,6) circle (1pt) node{};
	\path[->] (1.2,4) edge[bend left] (4-.2,4);
	\path[->] (1.2,4) edge[bend left] (8-.2,4);
	\path[->] (4.2,4) edge[bend left] (10-.2,4);
	\path[->] (4.2,4) edge[bend left] (12-.2,4);
	\path[->] (8.2,4) edge[bend left] (13-.2,4);
	\path[->] (8.2,4) edge[bend left] (15-.2,4);
	\path[->] (1.2,0) edge (5.8,0);
	\path[->] (6.2,0) edge (12.5-.2,0);
	\end{tikzpicture}
	\caption{Type $\mathscr{C}$: $\tilde w \to \tilde w \tilde w$. We begin with $C_1 \cup C_2$. In the first step, $C_1$ splits into $C_1 \cup C_3$. In the second step, $C_2$ splits into $C_2 \cup C_4$.}
	\label{ww-wwww}
\end{minipage} 
\efig

\bfig
\begin{minipage}{1.2\textwidth}
\centering

\includegraphics[height=2in]{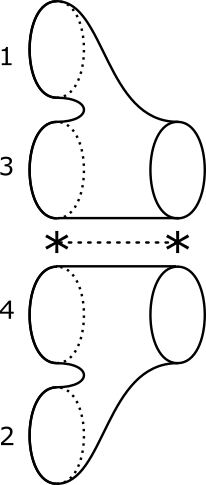}

\vspace{0.5in}

\hspace{-1in}
	\begin{tikzpicture}
	\filldraw 
	(0+22,2+4) circle (1pt) node[align=center, below=0.1 cm]{$w_+ w_-$}
	(1+22,0+4) circle (3pt) node[align=center, below=0.1 cm]{$w_- w_-$}
	(2+22,2+4) circle (1pt) node[align=center, below=0.1 cm]{$w_-w_+$}
	(1+22,4+4) circle (3pt) node[align=center, below=0.1 cm]{$w_+w_+$}
	(4+12,2+2) circle (1pt) node[align=center, below=0.1 cm]{$w_+w_-w_-$}
	(4+12,4+2) circle (1pt) node[align=center, below=0.1 cm]{$w_+w_+w_-$}
	(6+12,0+2) circle (1pt) node[align=center, below=0.1 cm]{$w_-w_-w_-$}
	(6+12,2+2) circle (1pt) node[align=center, below=0.1 cm]{$w_-w_+w_-$}
	(6+12,4+2) circle (1pt) node[align=center, below=0.1 cm]{$w_+w_-w_+$}
	(6+12,6+2) circle (1pt) node[align=center, below=0.1 cm]{$w_+w_+w_+$}
	(8+12,2+2) circle (1pt) node[align=center, below=0.1 cm]{$w_-w_-w_+$}
	(8+12,4+2) circle (1pt) node[align=center, below=0.1 cm]{$w_-w_+w_+$}
	(10-1,4)  circle (3pt) node[align=center, below=0.1 cm]{$w_+w_+ w_-w_-$}
	(12.5-1,0) circle (3pt) node[align=center, below=0.1 cm]{$w_- w_-w_- w_-$}
	(12.5-1,8) circle (3pt) node[align=center, below=0.1 cm]{$w_+w_+w_+w_+$}
	(15-1,4) circle (3pt) node[align=center, below=0.1 cm]{$w_- w_-w_+w_+$}
	(11-1,2) circle (1pt) node{}
	(11-1,4) circle (1pt) node{}
	(11-1,6) circle (1pt) node{}
	(12-1,2) circle (1pt) node{}
	(12-1,4) circle (1pt) node{}
	(12-1,6) circle (1pt) node{}
	(13-1,2) circle (1pt) node{}
	(13-1,4) circle (1pt) node{}
	(13-1,6) circle (1pt) node{}
	(14-1,2) circle (1pt) node{}
	(14-1,4) circle (1pt) node{}
	(14-1,6) circle (1pt) node{};
	\path[->] (12.5-1+0.2,8) edge (6+12-0.2,6+2);
	\path[->] (6+12+0.2,6+2) edge (1+22-0.2,4+4);
	\path[->] (10-1+0.2,4) edge[bend left] (6+12-0.2,2+2);
	\path[->] (6+12+0.2,2+2) edge[bend left]  (0+23-0.2,4);
	\path[->] (15-1+0.2,4) edge[bend left]  (8+12-0.2,2+2);
	\path[->] (8+12+0.2,2+2) edge[bend left]  (0+23-0.2,4);

	\end{tikzpicture}
	\caption{Type F: $\tilde w \tilde w \to \tilde w$. We start with $C_1 \cup C_2 \cup C_3 \cup C_4$, where $C_1 \cup C_2$ and $C_3 \cup C_4$ are each equivariant pairs of circles. In the first step, $C_1 \cup C_3$ merges to $C_1$ and $C_2$ stays put. In the second step, $C_2 \cup C_4$ merges to $C_2$ and $C_1$ stays put.}
	\label{wwww-ww}
\end{minipage} 

\efig


\epf
\eprop

\subsection{Visualizing higher differentials} \label{visualizing-higher-differentials}

Just as all components of $\hv d^2$ are induced by cancellation of some $\tau^\#$ arrow, longer differentials are induced by the cancellation of multiple $\tau^\#$ arrows. This can be seen by inducting on the applications of the cancellation operation that result in a differential of length $2r$. At each step, the slope of the arrows representing the induced differentials are each of the form $- \frac{k}{k-1}$. The generators at the ends of the staircase must be equivariant, and all the intermediate generators (i.e.\ the corners of the staircase) must be nonequivariant because they are adjacent to $\tau^\#$ arrows.

To prove that there are no higher differentials in the $\hvE$ spectral sequence, we  focus on individual staircases. This directs us to restrict our focus to a local piece of the complex surrounding the staircase which is large enough to determine the behavior of the map induced by the staircase.

\bdefn \label{defn-subcubes}

Let $\CD$ be an annular link diagram with $n$ crossings. 
Let $\alpha$ and $\beta$ be two $n$-bit strings where $\alpha<\beta$. 
The \emph{subcube of resolutions} from $R_\alpha$ to $R_\beta$, denoted $Cube(R_\alpha, R_\beta)$, is the subgraph of $Cube(\CD)$  whose 
\bi
\item vertex set consists of resolutions $R_\gamma$, where $\gamma$ is an $n$-bit string with
 $\gamma[i] = \alpha[i] = \beta[i]$ whenever $\alpha[i]= \beta[i]$
 \item edge set consists of all edges of $Cube(\CD)$ whose endpoints are both in the vertex set.
\ei

In $\scube(\CD)$, the \emph{subcube of chains} from $V_\alpha$ to $V_\beta$, denoted $\scube(V_\alpha, V_\beta)$, is defined similarly, with
\bi
\item underlying chain groups consisting of the vector spaces corresponding to the vertex set of $Cube(R_\alpha, R_\beta)$
\item differentials corresponding to the edge set of $Cube(R_\alpha, R_\beta)$.
\ei

In this case the crossings $c_i$ in the diagram $\CD$ for which $\alpha[i] \neq \beta[i]$ are called \emph{participating crossings}, and all other crossings are called \emph{nonparticipating crossings}. Let $\CD(\alpha, \beta)$ denote the link diagram obtained from $\CD(L)$ by smoothing only the nonparticipating crossings based on the common values of $\alpha$ and $\beta$ at these bits. Then $Cube(\CD(\alpha, \beta)) \cong Cube(R_\alpha, R_\beta)$ and $\scube(\CD(\alpha, \beta)) \cong \scube(V_\alpha, V_\beta)$.

Similarly, the components of $R_\alpha$ fall into two categories: 
the \emph{participating circles}, which are adjacent to the sites of participating crossings, and the  \emph{nonparticipating circles}, which are the components of $R_\alpha$ that do not participate in a merge or split along any edge in $Cube(R_\alpha, R_\beta)$.

\edefn

\bdefn \label{defn-equivariantly-split}
If $R_\alpha$ and $R_\beta$ are equivariant resolutions of an equivariant annular link diagram $\CD$, and $\alpha < \beta$, $Cube(R_\alpha, R_\beta)$ is \emph{equivariantly split} if the set of participating circles can be divided into two or more equivalence classes  under the equivalence relation 
\begin{quote}
$C_i \sim C_j$ if there is a (site of a) crossing adjacent to both $C_i$ and $C_j$. 
\end{quote}
\edefn

\bprop \label{no-long-diff-equiv-split} 
For an equivariantly split 2-periodic diagram $\CD$ with $2r$ crossings, there are no induced differentials of length $2r$.

\bpf 
If $\CD$ were equivariantly split, it could be written as a split diagram $\CD_1 \coprod \CD_2$, and hence
$${\text{\sc Cube}(\CD)} \cong \text{\sc Cube}(\CD_1) \otimes \text{\sc Cube}(\CD_2).$$ 
Let $\emptyset \neq \{C_j\}_{j \in J_1}$ denote the set of circles in ${R_{\alpha}}$ corresponding to the all-zeros resolution of $\CD_1$, and $\emptyset \neq \{C_j\}_{j \in J_2}$ the set of circles in ${R_{\alpha}}$ corresponding to the all-zeros resolution of $\CD_2$. 

The marked resolution $x$ corresponds to a choice of markings $a \in \{ v_-, v_+\}^{|J_1|}$ for $\{C_j\}_{j \in J_1}$ plus a choice of markings $\{ v_-, v_+\}^{|J_2|}$ for $\{C_j\}_{j \in J_2}$, so for computational reasons we write ${x} = a \otimes b$, even though in reality we might have ordered the tensor components differently. (We are implicitly using the isomorphism ${\text{\sc Cube}(\CD)} \cong \text{\sc Cube}(\CD_1) \otimes \text{\sc Cube}(\CD_2)$.)
Since ${x} = a \otimes b$ is an equivariant  generator, $a$ and $b$ must be equivariant  generators for $\text{\sc Cube}(\CD_1)$ and $\text{\sc Cube}(\CD_2)$, respectively.

As the differential respects the tensor factors, we may use the product rule to compute $d^{2r}({x})$. 
\begin{align*}
{x} &= a \otimes b \\
\partial {x} &= \partial a \otimes b + a \otimes \partial b \\
(1+\tau^\#)\inv \partial {x} &=
	(1+\tau^\#)\inv \partial a \otimes b + a \otimes (1+\tau^\#)\inv \partial b \\
\partial (1+\tau^\#)\inv \partial {x}  &=
	\partial (1+\tau^\#)\inv \partial a \otimes b + (1+\tau^\#)\inv \partial a \otimes \partial b + \partial a \otimes (1+\tau^\#)\inv \partial b \\
	& \indent + a \otimes \partial (1+\tau^\#)\inv \partial b\\   
\ldots
\end{align*}
Here we abuse notation and write $\partial$ to mean the differential in 
${\text{\sc Cube}(\CD)}$, $\text{\sc Cube}(\CD_1)$, or $ \text{\sc Cube}(\CD_2)$, depending on context.
Writing $\tilde \partial := (1+\tau^\#)\inv \partial$, where $(1+\tau^\#)\inv$ represents ``lifting by $1+\tau^\#$'' (see Remark \ref{cancellation-and-lift}), the above computation becomes
\begin{align*}
d^{2r} {x} &= \partial \tilde \partial^{2r-1} {x}  \\
	&= \partial \sum_{i = 0}^{2r-1} {{2r-1} \choose i }\tilde \partial^i a \otimes \tilde \partial^{2r-1-i} b \\
	&= \sum_{i = 0}^{2r-1} {{2r-1} \choose i } 
		\partial \tilde \partial^i a \otimes \tilde \partial^{2r-1-i} b 
		+ \tilde \partial^i a \otimes \partial \tilde \partial^{2r-1-i} b.\\
\end{align*}

We want to show that each summand must be 0 on $E^{2r}$, where all the survivors are necessarily equivariant.
Note that in order for $a' \otimes b' \in \text{\sc Cube}(\CD_1) \otimes \text{\sc Cube}(\CD_2)$ to be equivariant, $a'$ and $b'$ must be simultaneously equivariant, and hence must fall in a column corresponding to even Hamming weight in their respective {\sc Cubes}. 

\subsubsection*{Case $i$ odd} 
Then $\partial \tilde \partial^i a$ and $\tilde \partial^{2r-1-i}b$ have even Hamming weight,
while $\tilde \partial^i a$ and $\partial \tilde \partial^{2r -1 -i}b$ have odd Hamming weight. 
If any component of $\partial \tilde \partial^i a \otimes \tilde \partial^{2r-1-i}b$ is equivariant, then $\tilde \partial^{2r-1-i}b$ is equivariant, which contradicts the fact that $(1+\tau^\#) \tilde \partial^{2r-1-i}b = \partial^{2r -1 -i} \neq 0$. 

\subsubsection*{Case $i$ even}
Then $\partial \tilde \partial^i a$ and $\tilde \partial^{2r-1-i}b$ have odd Hamming weight,
while $\tilde \partial^i a$ and $\partial \tilde \partial^{2r -1 -i}b$ have even Hamming weight. 
If any component of  $\tilde \partial^i a \otimes \partial \tilde \partial^{2r -1 -i}b$ is equivariant, then $\tilde \partial^i a$ is equivariant, contradicting the fact that $(1+\tau ^\#)\tilde \partial^i a = \partial^i a \neq 0$. 

\epf

\eprop

\subsection{Grading shifts of higher differentials} \label{grading-shifts-of-higher-differentials}

Again let $\alpha, \beta$ be equivariant $2r$-bit strings with $\alpha < \beta$.
By Proposition \ref{no-long-diff-equiv-split}, we may now assume that ${\CD(\alpha, \beta)}$ is not equivariantly split. This is will allow us to perform computations by describing parts of the complex using connected graphs.

A first obstruction to the existence of a $^{hv}d^{2r}$ map between distinguished generators $x$ and $y$ is that $x$ and $y$ must have the same $j$-grading. For this reason, we compare the \emph{$j$-span} of ${V_{\alpha}}$ and ${V_{\beta}}$, defined as follows.

\bdefn  \label{defn-j-span}
The \emph{$j$-span} of the vector space ${V_\omega}$ associated with the binary string $\omega$ is the range of $j$-gradings of the distinguished generators of $V_{\omega}$. 
\edefn

Split maps shift the upper bound of the $j$-span up by 2, and merge maps shift the lower bound of the $j$-span up by 2. 

\bdefn \label{defn-path-through-cube}
A \emph{path} through a cube of resolutions is a sequence of vertices $R_{\alpha_i}$ such that $\alpha_{i+1}$ is a bit increment from $\alpha_i$, together with the edges between successive resolutions.
\edefn

Since any component of $^{hv}d^{2r}$ will come from repeatedly mapping via $\partial$ and lifting by $1+\tau^\#$ along a (directed) path from ${V_{\alpha}}$ to ${V_{\beta}}$ in $\text{\sc Cube}({\CD(\alpha, \beta)})$, 
and $1+\tau^\#$ affects neither of the $j$- or $k$-gradings, we focus on paths through $\scube({\CD(\alpha, \beta)})$ and which combinations of merges and splits they involve.

\blem \label{same-merge-split}
Any path through a cube $\sc{Cube}(\CD)$ from the all-zeros resolution to the all-ones resolution have the same number of merges and splits (separately).

\bpf
Let $\# m$ and $\# \Delta$ denote the number of merge and split maps along a path from ${V_{\alpha}}$ to ${V_{\beta}}$. These quantities satisfy the linear equations
$$ \# m + \# \Delta = 2r $$
$$ |R_{\beta}| - |R_{\alpha}| = \# \Delta - \# m. $$
\epf

\elem

So, in order to count merge and splits, we may without loss of generality choose any path we'd like. Choose $\gamma$ to be a path through the cube for which every other vertex is equivariant.
In other words, after every two steps along this path, we land on an equivariant resolution. 
Each such set of two-steps corresponds to one step in a path in $\scube({\CD(\alpha,\beta)}/\tau_A)$, since each equivariant pair of crossings corresponds to a crossing downstairs. We can thus classify the six different two-step paths in $\gamma$ by looking at the classification of edges in the cube downstairs.

\begin{figure} 
\centering 
\def\svgwidth{300pt} 
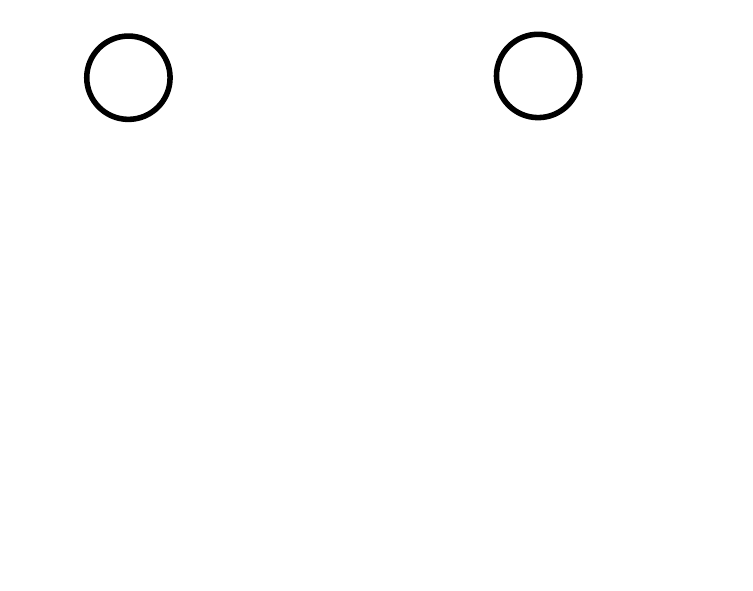
\caption{The six \emph{atoms} representing the six possible equivariant two-steps.}
\label{atoms}
\end{figure}

\bdefn \label{defn-atoms-molecules}
We can describe each two-step using modular pieces of a graph, which we'll call \emph{atoms} (Figure \ref{atoms}). These pieces are oriented in the sense that the lines on the left of the circle represent \emph{input circles}, or circles in the source resolution, and the lines on the right of the circle represent the \emph{output circles}, or circles in the target resolution of the two-step. 
A double bond indicates an equivariant twin pair of $W$ circles, and a single bond indicates a single equivariant $V$ circle.
While attaching these atoms together, we can only attach $V$ circles to $V$ circles, and equivariant pairs of $W$ circles to equivariant pairs of $W$ circles. 
A connected graph obtained by attaching these atoms together is called a \emph{molecule}.
\edefn

Since the ${\CD(\alpha, \beta)}/\tau_A$ diagram downstairs is nonsplit, when the pieces of $\gamma$ are attached at their equivariant ends, we obtain a single molecule $\Gamma$.
Let $A,B,C,D,E,F$ denote the number of modular pieces of type $\SA, \SB, \SC, \SD, \SE,$ and $\SF$ respectively.
In particular, this means that $A + B + C + D + E + F  = r$.

\bprop \label{eq-path-span}
Given an equivariant path $\gamma$ through an equivariantly nonsplit cube, the size of the $j$-span overlap between the vector spaces corresponding to the all-zeros resolution $V_\alpha$ and the all-ones resolution $V_\beta$ is 
$$ 2(A + B + 2C + D + E + 2F - N_b)$$
where
\bi
\item $A,B,C,D,E, $ and $F$ is the number of atoms of type $\SA, \SB, \SC, \SD, \SE,$ and $\SF$, respectively, in the molecule $\Gamma$ corresponding to $\gamma$
\item $N_b$ is the total number of bonds (counted with multiplicity) in $\Gamma$.
\ei

\bpf
Let $j_0$ denote the minimal $j$-grading of a distinguished generator of $V_{\alpha}$.
The number of components of ${R_{\alpha}}$ is equal to the number of input circles of the atoms separately, minus the number of input circles which are then identified with output circles of other atoms, which is equal to the number of total bonds, counted with multiplicity. Let $N_b$ denote the total number of bonds.
The source $j$-span is therefore 
$$ [ j_0, j_0 + 2( A+2B+2C+3D+2E+4F-N_b)].$$

Traveling through the path $\gamma$, we ultimately shift the bottom of the $j$-span by (two times) the number of merges, which is $2(B+2D+E+2F)$. 
Hence the possible $j$-span overlap between ${V_{\alpha}}$ and ${V_{\beta}}$ is 
$$
[j_0 + 2(B + 2D + E + 2F), j_0 + 2(A + 2B + 2C + 3D + 2E + 4F - N_b)]
$$
which is a difference of 
$2(A + B + 2C + D+ E+ 2F - N_b)$.

\epf

\eprop

Thus the largest $j$-span overlap using the atoms in $\Gamma$ is achieved when $\Gamma$ has a minimal total number of bonds.

\blem \label{tree-minimizes-bonds}
The minimal total number of bonds achievable by a connected molecule $\Gamma$ comprised of a fixed set of atoms is achieved when $\Gamma$ is composed to two trees, one containing all the $\SA,\SB,\SD,$ and $\SE$ atoms, and the other containing all the $\SC$ and $\SF$ atoms, with the two trees connected by a unique double bond.

\bpf

The structure of a minimally-bonded $\Gamma$ must be a tree: if $\Gamma$ has a cycle, we can decrease the number of bonds by cutting any side of the cycle.

Furthermore, when attaching atoms of type $\SA,\SB,\SD,$ and $\SE$, the number of bonds is minimized by connecting to the molecule via a single rather than a double bond. When these atoms are arranged as a subtree of $\Gamma$, the number of single bonds is minimized. 
The $\SC$ and $\SF$ atoms will invariably attach to the tree via a double bond, so the two-tree molecule described has the same number of double bonds as any minimally-bonded $\Gamma$.
\epf

\elem

\subsection{Higher differentials vanish} \label{higher-differentials vanish}

The final major step in the proof of Theorem \ref{main-thm} uses the $j$- and $k$- gradings to show that the $\hvE$ spectral sequence collapses on $\hvE^3$.

\bprop \label{higher-differentials-vanish}
$^{hv} d^{2r} = 0$ for all $r \geq 2$.
\bpf

By Proposition \ref{no-long-diff-equiv-split}, we may assume that $Cube(R_\alpha, R_\beta)$ is equivariantly nonsplit. Pick an equivariant path $\gamma$ through $Cube(R_\alpha, R_\beta)$, and let $\Gamma$ denote the associated molecule.
Let's first assume that at least one of the atoms of $\Gamma$ is \emph{not} of type $\SC$ or $\SF$. 
In the worst case, the total number of bonds is the sum of 
$A + B+D+E-1$ from the $\mathscr{A,B,D,E}$ tree, $2C +2F-2$ from the $\mathscr{C,F}$ tree, and $2$ from connecting the two trees. Hence
$N_b = A + B + 2C + D + E + 2F - 1$.
The size of the overlap is then 
$$
 2(A + 2B + 2C + 3D + 2 E + 4F - (A + B + 2C + D + E + 2F - 1)) - 2(B+2D+E + 2F)
 = 2.$$ 
Therefore ${x}$ must have either the maximum or next-to-maximum $j$-grading of the distinguished generators of ${V_{\alpha}}$.
Since ${x}$ is equivariant, any equivariant pair of trivial circles must have the same marking, so since ${x}$ can have at most one component marked with a $v_-$, all trivial circles in ${x}$ are marked $w_+$. 
Moreover, from Figures \ref{v-vww} through \ref{wwww-ww} we see that if any equivariant paths out of ${V_{\alpha}}$ are not of type $\SE$, then $x$ is not in $\ker(\dtwo)$. So all equivariant paths out of ${V_{\alpha}}$ must be of type $\SE$. This indicates that ${x}$ is a marked resolution consisting solely of $v_+$-marked concentric circles. 
But this means ${x}$ is unique in its $k$-grading in the entire $\scube(V_\alpha, V_\beta)$. Therefore there cannot be a $\hv d^{2r}$ component from ${x}$ to ${y}$ in $\scube(V_\alpha, V_\beta)$.

In the case where \emph{all} the atoms are of type $\SC$ or $\SF$, the worst case $\Gamma$ has $2C + 2F - 2$ bonds, producing a $j$-span overlap of 
$$
2(2C + 4F - (2C + 2F - 2)) - 2(2F) = 4.$$
Therefore ${x}$ must have at most two circles marked with $w_-$. Again, we see in Figures \ref{ww-wwww} and \ref{wwww-ww} that ${x}$ is not in $\ker(\dtwo)$.

Hence there cannot be a $^{hv}d^{2r}$ component from (a $\theta$-power translate of) ${x}$ to (a $\theta$-power translate of) ${y}$.
\epf

\eprop

As $\atate(\CD(\tilde L))$ is 1-periodic in the horizontal direction, its homology is 1-periodic as well, so the total rank along each diagonal corresponding to the induced filtration gradings on $H^{Tot}$ is the same as the total rank down along a vertical column. 
Therefore $AKh(L) \otimes_{\FF} \F[\theta, \theta\inv] 
\cong \hvE^3 = \hvE^\infty \cong \vhE^\infty$ so the $\vhE$ spectral sequence has $\vhE^1 \cong AKh(\tilde L) \otimes_{\FF} \F[\theta, \theta\inv]$ and $\vhE^\infty \cong AKh(L) \otimes_{\FF} \F[\theta, \theta\inv]$. 

Since all the maps in the double complex preserve $j$- and $k$-grading, the spectral sequence splits along $j$- and $k$-gradings. This concludes the proof of Theorem \ref{main-thm}. 

Finally, since each page of a spectral sequence is a subquotient of the previous page, Corollary \ref{main-cor} follows:

\begin{repcorollary}{main-cor}

For any annular tangle closure $L$, quantum grading $j$ and $\sltwo$ weight space  grading $k$,
$$
\rk AKh^{j,k}(L) \leq \rk AKh_{2j-k, k} (\tilde L).
$$

\end{repcorollary}

\section{Example: Hopf link and stabilized unknot} \label{eg-hopf-link}

To help convey the construction to the reader, we describe the simplest nontrivial example, a two-crossing annular Hopf link as the 2-braid closure $\tilde L :=\widehat{\sigma_1^2}$. The quotient link is a positively stabilized unknot $L = \widehat{\sigma_1}$.

Annular diagrams for both $\tilde L$ and $L$ are shown in Figure \ref{annular-hopf}.
We label the crossings of $\CD(\tilde L)$ in accordance to our convention (\S \ref{involution-on-ckh}); the first crossing is on the right, and the second crossing is on the left. The cube of resolutions $Cube(\CD(\tilde L))$ is shown in Figure \ref{hopf-res-cube}. In order to simplify the notation in this small example, we give the distinguished generators unique, concise names in the cube of chains shown in Figure \ref{hopf-vs-cube-2}. For example, ``$a_{+-}$'' is the generator $v_+ \otimes v_-$ in $V_{00}$, corresponding to the marked resolution shown in Figure \ref{a-pm-marked-res}.

\bfig
\centering
\def\svgwidth{400pt} 
\begingroup%
  \makeatletter%
  \providecommand\color[2][]{%
    \errmessage{(Inkscape) Color is used for the text in Inkscape, but the package 'color.sty' is not loaded}%
    \renewcommand\color[2][]{}%
  }%
  \providecommand\transparent[1]{%
    \errmessage{(Inkscape) Transparency is used (non-zero) for the text in Inkscape, but the package 'transparent.sty' is not loaded}%
    \renewcommand\transparent[1]{}%
  }%
  \providecommand\rotatebox[2]{#2}%
  \ifx\svgwidth\undefined%
    \setlength{\unitlength}{349.15188767bp}%
    \ifx\svgscale\undefined%
      \relax%
    \else%
      \setlength{\unitlength}{\unitlength * \real{\svgscale}}%
    \fi%
  \else%
    \setlength{\unitlength}{\svgwidth}%
  \fi%
  \global\let\svgwidth\undefined%
  \global\let\svgscale\undefined%
  \makeatother%
  \begin{picture}(1,0.69784902)%
    \put(-0.01433837,1.16624286){\color[rgb]{0,0,0}\makebox(0,0)[lt]{\begin{minipage}{0.16038865\unitlength}\raggedright \end{minipage}}}%
    \put(0,0){\includegraphics[width=\unitlength,page=1]{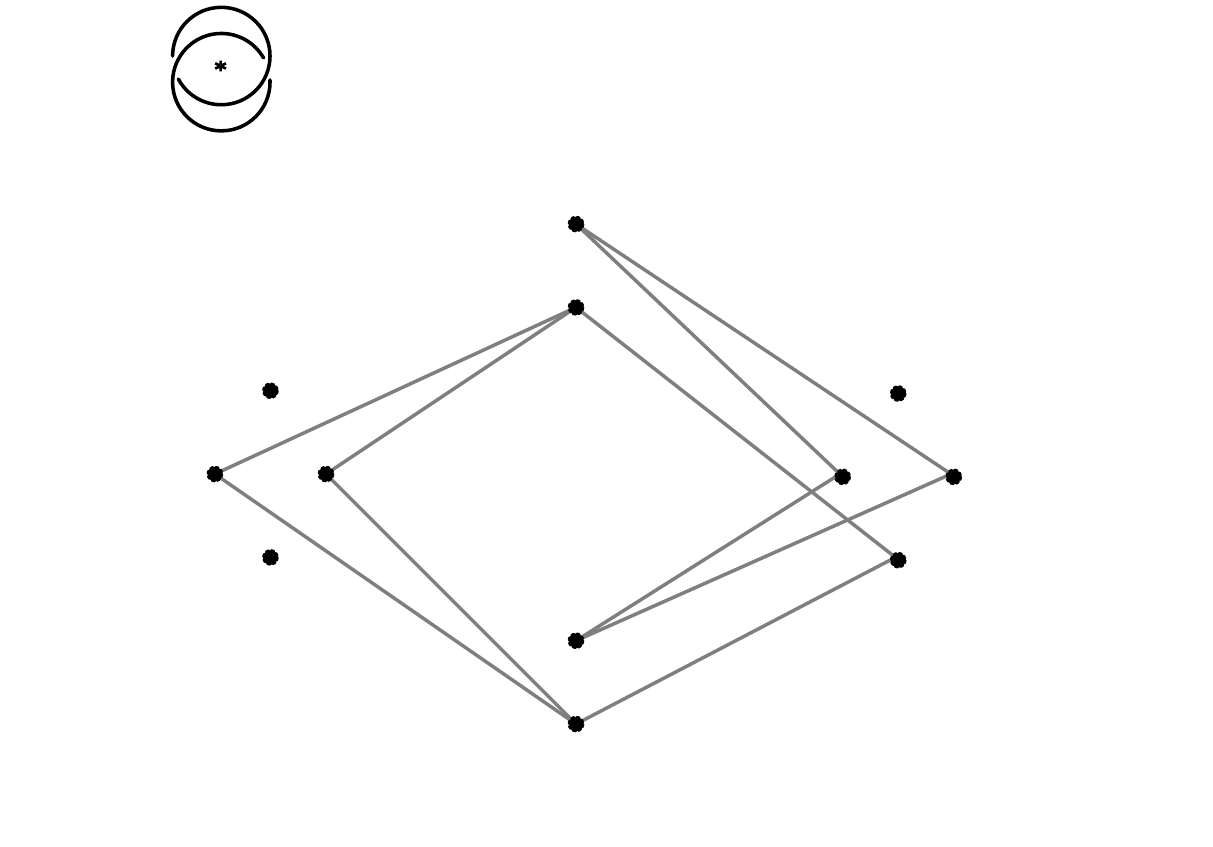}}%
    \put(0.18641258,0.34578427){\color[rgb]{0,0,0}\makebox(0,0)[lb]{\smash{$a_{++}$}}}%
    \put(0.25285929,0.26788121){\color[rgb]{0,0,0}\makebox(0,0)[lb]{\smash{$a_{-+}$}}}%
    \put(0.12340275,0.27246375){\color[rgb]{0,0,0}\makebox(0,0)[lb]{\smash{$a_{+-}$}}}%
    \put(0.20016018,0.21289082){\color[rgb]{0,0,0}\makebox(0,0)[lb]{\smash{$a_{--}$}}}%
    \put(0.45678202,0.47867772){\color[rgb]{0,0,0}\makebox(0,0)[lb]{\smash{$b_+$}}}%
    \put(0.46136455,0.3984834){\color[rgb]{0,0,0}\makebox(0,0)[lb]{\smash{$b_-$}}}%
    \put(0.46136455,0.14644409){\color[rgb]{0,0,0}\makebox(0,0)[lb]{\smash{$c_+$}}}%
    \put(0.46136455,0.06854103){\color[rgb]{0,0,0}\makebox(0,0)[lb]{\smash{$c_-$}}}%
    \put(0.73860772,0.343493){\color[rgb]{0,0,0}\makebox(0,0)[lb]{\smash{$d_{++}$}}}%
    \put(0.77297672,0.27017248){\color[rgb]{0,0,0}\makebox(0,0)[lb]{\smash{$d_{-+}$}}}%
    \put(0.64352016,0.27590064){\color[rgb]{0,0,0}\makebox(0,0)[lb]{\smash{$d_{+-}$}}}%
    \put(0.71569505,0.19685195){\color[rgb]{0,0,0}\makebox(0,0)[lb]{\smash{$d_{--}$}}}%
    \put(0.71569505,0.20372579){\color[rgb]{0,0,0}\makebox(0,0)[lb]{\smash{}}}%
    \put(0.02315995,0.62761008){\color[rgb]{0,0,0}\makebox(0,0)[lb]{\smash{$\scube($}}}%
    \put(0.26431574,0.62761008){\color[rgb]{0,0,0}\makebox(0,0)[lb]{\smash{$)=$}}}%
  \end{picture}%
\endgroup%

\caption{The cube of chains for the annular Hopf link from Figure \ref{annular-hopf}, with $AKh$ differentials drawn as solid lines and $\tau^\#$ drawn as dotted lines. The arrow heads have been dropped; $gr_i$ increases to the right. This is a finer version of Figure \ref{hopf-vs-cube}.}
\label{hopf-vs-cube-2}
\efig

\bfig
\centering
\def\svgwidth{100pt} 
\begingroup%
  \makeatletter%
  \providecommand\color[2][]{%
    \errmessage{(Inkscape) Color is used for the text in Inkscape, but the package 'color.sty' is not loaded}%
    \renewcommand\color[2][]{}%
  }%
  \providecommand\transparent[1]{%
    \errmessage{(Inkscape) Transparency is used (non-zero) for the text in Inkscape, but the package 'transparent.sty' is not loaded}%
    \renewcommand\transparent[1]{}%
  }%
  \providecommand\rotatebox[2]{#2}%
  \ifx\svgwidth\undefined%
    \setlength{\unitlength}{58.05983536bp}%
    \ifx\svgscale\undefined%
      \relax%
    \else%
      \setlength{\unitlength}{\unitlength * \real{\svgscale}}%
    \fi%
  \else%
    \setlength{\unitlength}{\svgwidth}%
  \fi%
  \global\let\svgwidth\undefined%
  \global\let\svgscale\undefined%
  \makeatother%
  \begin{picture}(1,0.98863528)%
    \put(0,0){\includegraphics[width=\unitlength,page=1]{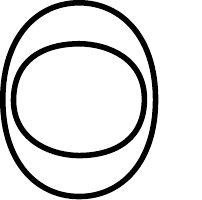}}%
    \put(0.35967658,0.40692407){\color[rgb]{0,0,0}\makebox(0,0)[lb]{\smash{*}}}%
    \put(0.72146099,0.85559949){\color[rgb]{0,0,0}\makebox(0,0)[lb]{\smash{$+$}}}%
    \put(0.5044436,0.56968755){\color[rgb]{0,0,0}\makebox(0,0)[lb]{\smash{$-$}}}%
    \put(-0.10527151,4.3003213){\color[rgb]{0,0,0}\makebox(0,0)[lt]{\begin{minipage}{1.24009999\unitlength}\raggedright \end{minipage}}}%
    \put(-1.55205478,2.85353801){\color[rgb]{0,0,0}\makebox(0,0)[lt]{\begin{minipage}{5.99381651\unitlength}\raggedright \end{minipage}}}%
  \end{picture}%
\endgroup%

\caption{The marked resolution corresponding to the generator $a_{+-}$.}
\label{a-pm-marked-res}
\efig

Putting this all together, the bicomplex $\atate(\CD(\tilde L))$ locally looks like Figure \ref{tate-cpx-eg}.

\bfig
\def\svgwidth{400pt} 
\begingroup%
  \makeatletter%
  \providecommand\color[2][]{%
    \errmessage{(Inkscape) Color is used for the text in Inkscape, but the package 'color.sty' is not loaded}%
    \renewcommand\color[2][]{}%
  }%
  \providecommand\transparent[1]{%
    \errmessage{(Inkscape) Transparency is used (non-zero) for the text in Inkscape, but the package 'transparent.sty' is not loaded}%
    \renewcommand\transparent[1]{}%
  }%
  \providecommand\rotatebox[2]{#2}%
  \ifx\svgwidth\undefined%
    \setlength{\unitlength}{324.26126862bp}%
    \ifx\svgscale\undefined%
      \relax%
    \else%
      \setlength{\unitlength}{\unitlength * \real{\svgscale}}%
    \fi%
  \else%
    \setlength{\unitlength}{\svgwidth}%
  \fi%
  \global\let\svgwidth\undefined%
  \global\let\svgscale\undefined%
  \makeatother%
  \begin{picture}(1,0.72711131)%
    \put(0,0){\includegraphics[width=\unitlength,page=1]{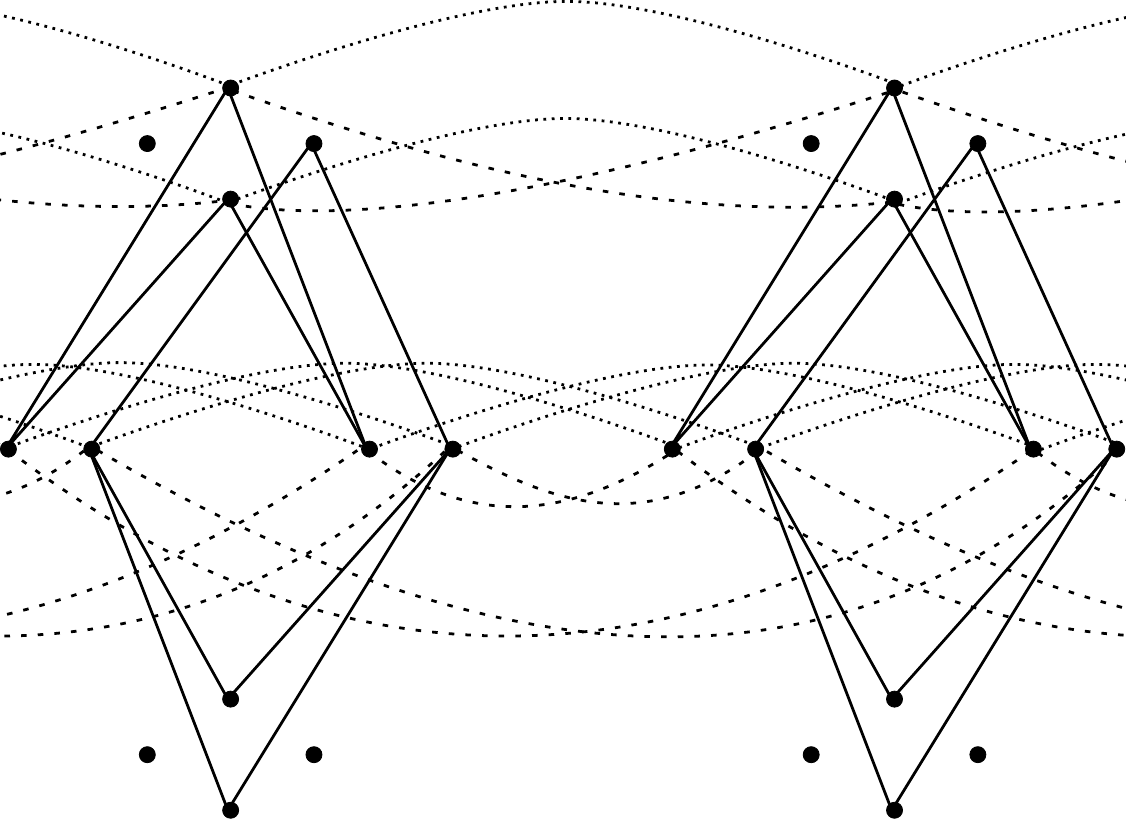}}%
  \end{picture}%
\endgroup%

\caption{A local picture of the $\atate(\CD(\tilde L))$ bicomplex for the Hopf link example. There are two copies of the cube of chains in this figure, rotated $90^\circ$ clockwise from from Figure \ref{hopf-vs-cube}. The dotted lines which bend up represent the ``identity'' map which increases the $t$-grading, and the dashed lines which bend down represent the $\tau^\#$ map. This local picture is repeated horizontally to form $\atate(\CD(\tilde L))$.}
\label{tate-cpx-eg}
\efig

After cancelling all the $\tau^\#$ arrows, the only remaining interesting page to compute the homology of is $\hvE^2$, shown in Figure \ref{tate-cpx-eg-d2}. As expected, each slope-(-2) line in $\hvE^2$ corresponds to the cube of chains for the quotient link $L$, under the correspondence described by Propositions \ref{E2-correspondence} and \ref{d2-correspondence}.

\bfig
\def\svgwidth{400pt} 
\begingroup%
  \makeatletter%
  \providecommand\color[2][]{%
    \errmessage{(Inkscape) Color is used for the text in Inkscape, but the package 'color.sty' is not loaded}%
    \renewcommand\color[2][]{}%
  }%
  \providecommand\transparent[1]{%
    \errmessage{(Inkscape) Transparency is used (non-zero) for the text in Inkscape, but the package 'transparent.sty' is not loaded}%
    \renewcommand\transparent[1]{}%
  }%
  \providecommand\rotatebox[2]{#2}%
  \ifx\svgwidth\undefined%
    \setlength{\unitlength}{324.26126862bp}%
    \ifx\svgscale\undefined%
      \relax%
    \else%
      \setlength{\unitlength}{\unitlength * \real{\svgscale}}%
    \fi%
  \else%
    \setlength{\unitlength}{\svgwidth}%
  \fi%
  \global\let\svgwidth\undefined%
  \global\let\svgscale\undefined%
  \makeatother%
  \begin{picture}(1,0.72711131)%
    \put(0,0){\includegraphics[width=\unitlength,page=1]{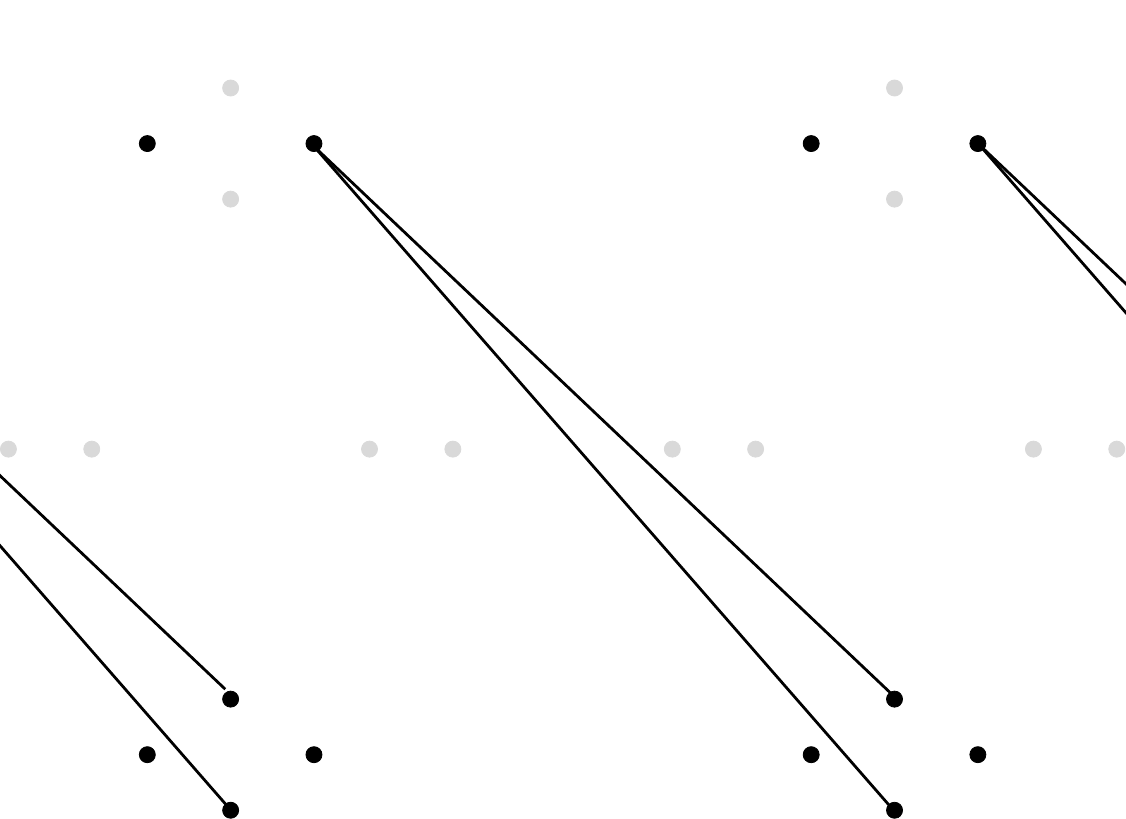}}%
  \end{picture}%
\endgroup%

\caption{A local picture of the page $\hvE^2$ in our Hopf link example.}
\label{tate-cpx-eg-d2}
\efig

In this case, the fact that there are no higher differentials is clear from the fact that $CKh(\CD(\tilde L))$ is only supported on three adjacent $i$-gradings. For comparison, the annular Khovanov homologies for $\tilde L$ and $L$ are listed in the following chart.

\begin{center}
\begin{tabular}{|c | c | c | c|}
\hline 
$(j,k)$ & $\rk (AKh^{j,k}(L)) $ & $(2j-k, k)$ & $\rk(AKh_{2j-k, k}(\tilde L))$ \\
\hline
(3,2) & 1 & (4,2) & 1 \\
\hline
(1,0) & 1 & (2,0) & 1 \\
\hline
(-1, -2) & 1 & (0,-2) & 1 \\
\hline
(3,0) & 1 & (6,0) & 1 \\
\hline
& & (4,0) & 2 \\
\hline
\end{tabular}
\end{center}

\section{Decategorification} \label{decat}

Annular Khovanov homology was originally described as a categorification of the Kauffman bracket skein module in the trivial $I$-bundle over the annulus.
In our annular setting, the skein module is an $\FF[q^{\pm 1}]$-module
consisting of all links in $A \times I$ modulo the usual Kauffman bracket skein relations.
The set of all trivial braid closures $\CB = \{ \widehat{\mathbbm{1}}_s \st s =0, 1, 2, \ldots\}$ form a canonical basis for the skein module. 

Following the conventions in \cite{Roberts-2013}, we use the variables $t,q,$ and $x$ to record the $i$-, $j$-, and $k$-gradings, respectively. Setting $q = A^{-2}$ and introducing the $k$-grading to the skein relation, we have a variant of the Kauffman bracket skein relations. Choose a particular crossing $c$ and let $L_0$ and $L_1$ denote the links obtained from $L$ by smoothing $c$ by the 0-resolution and 1-resolution, respectively. Then our relations are 
\begin{align*}
\langle L \rangle & = \langle L_0 \rangle + tq  \langle L_1 \rangle \\
\langle L \cup U \rangle & = (q + q\inv) \langle L \rangle
\end{align*}
where $U$ is an annular unknot (a trivial circle). 

Normalizing by setting $\widehat{\mathbbm{1}}_1 \mapsto qx + q\inv x\inv$, we obtain a map 
$$ \phi: \FF[t, q^{\pm 1}, \CB] \to \FF[t,q^{\pm1},x^{\pm1}]$$
by viewing $\CB$ as a monoid under disjoint union. 
Finally, we add in the proper grading shifts to obtain a variant of the Jones polynomial over $\FF$ by setting $V(t,q,x) = t^{n_-}q^{n_+ - 2n_-} \phi(L)$. (To be clear, one needs to first orient $L$ to obtain $n_-$ and $n_+$, and we are viewing the link $L$ as an element of the skein module.) 

Of interest to us is the annular link invariant $V_L(-1,q,x)$ which is the graded Euler characteristic of $AKh(L)$. (See Theorem 8.1 in \cite{APS-2004} for details.) 
In other words, 
$$ 
V_L(t,q,x) = \sum_{i,j,k} t^i q^j x^k \rk (AKh^{i,j,k}(L))
$$
so the graded Euler characteristic is 
$$
V_L(-1,q,x) 
= \sum_{i,j,k} (-1)^i q^j x^k \rk (AKh^{i,j,k}(L)).$$
As in \cite{APS-2004} and \cite{Roberts-2013}, we view the Kauffman bracket as a set of polynomials by defining $q_{k,L}(t,q)$ to be the coefficient of $x^k$ in $V(t,q,x)$, so that 
$$q_{k,L}(-1,q) = \sum_{i,j} (-1)^i q^j \rk AKh^{i,j,k}(L)$$
 is the coefficient of $x^k$ in the decategorification of $AKh(L)$.\\

The decategorification of Theorem \ref{main-thm} is the following:

\begin{repcorollary}{decat-cor}
For all $j$ and $k$, 
$$\langle q_{k, \tilde L}(-1,q), q^{2j-k} \rangle \equiv \langle q_{k,L}(-1,q), q^j \rangle \mod 2$$
where $\langle f,g \rangle$ denotes the coefficient of $g$ in $f$. 

\bpf
This follows from Theorem \ref{main-thm} and the fact that cancellation reduces rank by a multiple of 2 at each step (represented by erasing the dots at the head and tail of the arrow). 
\epf
\end{repcorollary}

In \cite{Murasugi-1971}, Murasugi proved a relationship between the Alexander polynomials of a periodic link and that of its quotient; in \cite{Murasugi-1988}, he proved an analogous formula for the Jones polynomial. By considering Theorem \ref{main-thm} with respect to the $j_1 = j-k$ grading rather than the $j$- and $k$-gradings individually,  we obtain an annular analogue to Murasugi's formulas: 

\begin{repcorollary}{murasugi-like}
$V_{\tilde L}(1,q, q\inv) \equiv [V_L(1,q,q\inv)]^2 \mod 2.$
\end{repcorollary}

If on the other hand one is more interested in the usual Jones polynomial of the 2-periodic link, we may forget the $k$-grading information by setting $x=1$, obtaining the following relationship:

\begin{repcorollary}{decat-kh}
$V_{\tilde L}(1,q,1) \equiv V_L(1,q^2,q\inv) \mod 2.$
\end{repcorollary}

It is then natural to ask if there is a categorification of this relationship, coming from a spectral sequence from the Khovanov homology of $\tilde L$ to the annular Khovanov homology of $L$ (with some grading information sacrificed). In the next section, we show that this spectral sequence indeed exists for some families of links, and is likely to exist in general by way of a bicomplex very similar to the $\atate$ bicomplex.

\section{The Khovanov-Tate bicomplex} \label{section-conjecture}

We can similarly define the \emph{Khovanov-Tate bicomplex} $\khtate(\CD(\tilde L))$ by replacing the vertical differentials $\partial_{AKh}$ in $\atate(\CD(\tilde L))$ with the Khovanov differential $\partial_{Kh}$.

\begin{repconjecture}{main-conj}
Let $\tilde L$ be a 2-periodic link in $S^3$ with quotient link $L$. There is a spectral sequence with
$$
E^1 \cong Kh(\tilde L) \otimes_{\FF} \FF[\theta, \theta\inv]
\abuts
E^\infty \cong AKh(L) \otimes_{\FF} \FF[\theta, \theta\inv].
$$
This would in turn imply the following cascade of rank inequalities:
$$
\rk AKh(\tilde L) 
\geq \rk Kh(\tilde L)
\geq \rk AKh(L)
\geq \rk Kh(L)
$$
where the first and third inequalities are given by the $k$-grading filtration on $CKh(\CD(\tilde L))$ and $CKh(\CD(L))$.
\end{repconjecture}

In \S \ref{khtate-analogous} we give abbreviated proofs of the $Kh$ analogues of the results from the $AKh$ case which support Conjecture \ref{main-conj}. But first, we discuss holistic reasons why we would expect this conjecture to hold.

\brmk \label{why-expect-this}
Aside from the proofs of theorems given in \S \ref{khtate-analogous},  this conjecture would not be hard to believe. 
$H_*(Tot(\khtate))$ is supposed to capture information about both the periodic link and the involution witnessing the symmetry, so the presence of the braid axis in $S^3$ is inherently captured in the data of the $\khtate$ bicomplex. 
Hence $H_*(Tot(\khtate))$ \emph{should} take on the form of an annular version of Khovanov homology.
Furthermore, Seidel and Smith proved the rank inequality for symplectic Khovanov homology \cite{Seidel-Smith-2010}, which is known to be equivalent to Khovanov homology over ground rings of characteristic 0 \cite{Abouz-Smith-2015}: 
$$ \rk Kh_{symp}(\tilde L) \geq \rk Kh_{symp}(L).$$ 
They conjecture that symplectic Khovanov homology is also isomorphic to Khovanov homology in characteristic 2 \cite{Seidel-Smith-2006}, so the rank inequalities implied by this conjecture agree with their result.
\ermk

\subsection{Analogous statements in the $\khtate$ case} \label{khtate-analogous}

Khovanov homology can be viewed as a deformation of annular Khovanov homology, defined for knots in $S^3$. This is why we have been denoting the underlying vector space for the annular Khovanov complex of a link diagram by $CKh$; this is the same underlying vector space as in the Khovanov complex.  Since there is no basepoint in the link diagram, the Khovanov distinguished generators are only doubly graded, by the homological $i$-grading and the quantum $j$-grading. There are only two types of edge maps: merge and split. Their definitions are shown in Figure \ref{kh-maps} using the dot and arrow notation we used to describe the $AKh$ differentials in back in \S \ref{annular-khovanov-homology}. 

\begin{figure}
\centering

\begin{minipage}{0.47\textwidth}
	\centering
	\begin{tikzpicture}
	\filldraw 
	(0,0) circle (2pt) node[align=center, below = 0.1 cm]{$x_-$} 
	(0,2) circle (2pt) node[align = center, below = 0.1 cm]{$x_+$}
	(4,4) circle (2pt) node[align=center, below = 0.1 cm]{$x_+\otimes x_+$}
	(3,2) circle (2pt) node[align=center, below = 0.1 cm]{$x_+ \otimes x_-$} 
	(5,2) circle (2pt) node[align=center, below = 0.1 cm]{$x_- \otimes x_+$} 
	(4,0) circle (2pt) node[align=center, below = 0.1 cm]{$x_- \otimes x_-$};
	\path[->] (0.2,0) edge (3.8,0);
	\path[->] (0.2,2) edge[bend left] (2.8, 2);
	\path[->] (0.2,2) edge[bend left] (4.8,2);
	\end{tikzpicture}
	\caption*{Merge map $m$}
	\label{kh-merge}
\end{minipage} 
\begin{minipage}{0.47\textwidth}
	\centering
\begin{tikzpicture}
\filldraw 
(8,2) circle (2pt) node[align=center, below = 0.1 cm]{$x_-$} 
(8,4) circle (2pt) node[align = center, below = 0.1 cm]{$x_+$}
(4,4) circle (2pt) node[align=center, below = 0.1 cm]{$x_+\otimes x_+$}
(3,2) circle (2pt) node[align=center, below = 0.1 cm]{$x_+ \otimes x_-$} 
(5,2) circle (2pt) node[align=center, below = 0.1 cm]{$x_- \otimes x_+$} 
(4,0) circle (2pt) node[align=center, below = 0.1 cm]{$x_- \otimes x_-$};
\path[->] (3.2,2) edge[bend left] (7.8, 2);
\path[->] (5.2,2) edge[bend left] (7.8, 2);
\path[->] (4.2, 4) edge (7.8, 4);
\end{tikzpicture} 
	\caption*{Split map $\Delta$}
	\label{kh-split}
\end{minipage} 

\caption{The Khovanov differentials. These do not depend on the basepoint $*$, so $x$ can stand for $v$ or $w$.}
\label{kh-maps}
\end{figure}

Let $\hvE$ and $\vhE$ denote the spectral sequences induced by the row-wise and column-wise filtrations, respectively, of the bicomplex $\khtate(\CD(\tilde L)) = CKh(\CD(\tilde L)) \otimes_{\FF} \FF[\theta, \theta\inv]$, where $\partial^v = \partial_{Kh}$ and $\partial^h = 1 + \tau^\#$.

\blem
$\partial_{Kh}$ commutes with $\tau^\#$.
\bpf
We need to check that the two types of edge maps (merge and split) commute with $\tau^\#$. This is easily checked, and clear from thinking about the involution acting on the saddle cobordism representing the merge or split map, just as in the $AKh$ cases.

\epf
\elem

\blem 
$\vhE^1 \cong Kh(\tilde L) \otimes_{\FF} \FF[\theta, \theta\inv]$.
\bpf
This follows from the same argument as in the proof of Lemma \ref{vhE1}.
\epf
\elem

\blem \label{kh-hvE1-gen-by-equiv-gens}
$\hvE^1$ is generated by the equivariant generators of $\khtate(\CD(\tilde L))$. Hence if $r$ is odd, $^{hv}d^r = 0$.
\bpf
Since the rows of the $\khtate(\CD(\tilde L))$ bicomplex are identical to those of the $\atate(\CD(\tilde L))$ complex, this is equivalent to Lemma \ref{hvE1-gen-by-equiv-gens}.
The second statement follows from the first just as in the proof of Lemma \ref{dr-odd-zero}.
\epf
\elem

\blem 
Every component of $^{hv} d^2$ is induced by the cancellation of a $\tau^\#$ arrow at a row corresponding to an odd Hamming weight.
\bpf
The proof of Lemma \ref{d2-from-cancellation} also holds in this case.
\epf
\elem

\bprop \label{kh-E2-correspondence}
There is a one-to-one correspondence

\begin{align*}
\left \{ \text{generators of } \hvE^2 \right \}  \longleftrightarrow &
\left \{ \text{generators of } CKh(L) \otimes_{\FF} \FF[\theta, \theta\inv] \right \}\\
\tilde v  \longleftrightarrow &  v
\end{align*}

induced by $\tau$.
\bpf
By Lemma \ref{kh-hvE1-gen-by-equiv-gens}, we have already identified the generators of $\hvE^2$ with the equivariant marked resolutions, which the proof of Proposition \ref{E2-correspondence} shows to be in bijection with the generators of $CKh(\CD(L)) \otimes_{\FF} \FF[\theta, \theta\inv]$ (since this is the underlying vector space in both   the $\khtate(\CD(L))$ and $\atate(\CD(L))$ bicomplexes).
\epf
\eprop

\bprop \label{kh-d2-correspondence}
Under the correspondence above, each slope -2 line in the bicomplex $(\hvE^2, ^{hv} d^2)$ is isomorphic to the complex $(CKh(\CD(L)), \partial_{Kh})$. 

\bpf
Again, we just need to compute the induced length 2 differentials. Since there are six types of $AKh$ differentials, we can verify them case-by-case. The computations are shown in Figures \ref{kh-v-vww} through \ref{kh-wwww-ww}.


\bfig
\begin{minipage}{\textwidth}
\centering
\includegraphics[height=2in]{v-vww.png}
	\begin{tikzpicture}
	\filldraw 
	(0,0) circle (3pt) node[align=center, below=0.1 cm]{$v_-$}
	(0,2) circle (3pt) node[align=center, below=0.1 cm]{$v_+$}
	(2,2) circle (1pt) node[align=center, below=0.1 cm]{$v_+ w_-$}
	(3,0) circle (1pt) node[align=center, below=0.1 cm]{$v_-  w_-$}
	(3,4) circle (1pt) node[align=center, below=0.1 cm]{$v_+  w_+$}
	(4,2) circle (1pt) node[align=center, below=0.1 cm]{$v_- w_+$}
	(6,2) circle (3pt) node[align=center, below=0.1 cm]{$v_+ w_-  w_-$}
	(6,4) circle (1pt) node[align=center, below=0.1 cm]{$v_+  w_+  w_-$}
	(8,0) circle (3pt) node[align=center, below=0.1 cm]{$v_-  w_-  w_-$}
	(8,2) circle (1pt) node[align=center, below=0.1 cm]{$v_-  w_+  w_-$}
	(8,4) circle (1pt) node[align=center, below=0.1 cm]{$v_+  w_-  w_+$}
	(8,6) circle (3pt) node[align=center, below=0.1 cm]{$v_+  w_+  w_+$}
	(10,2) circle (1pt) node[align=center, below=0.1 cm]{$v_- w_-  w_+$}
	(10,4) circle (3pt) node[align=center, below=0.1 cm]{$v_-  w_+  w_+$};
	\path[->] (0.2,2) edge[bend left] (1.8,2);
	\path[->] (0.2,2) edge[bend left] (3.8, 2);
	\path[->] (0.2,0) edge (2.8,0);
	\path[->] (2.2,2) edge[bend left] (5.8,2);
	\path[->] (4.2,2) edge[bend left] (7.8,2);
	\path[->] (4.2,2) edge[bend left] (9.8,2);
	\path[->] (3.2,0) edge (7.8,0);
	\end{tikzpicture}
	\caption{Type $\mathscr{A}$: $\tilde v \to \tilde v \tilde w$ for Khovanov differentials.}
	\label{kh-v-vww}
\end{minipage} 

\vspace{0.5 in}
\begin{minipage}{\textwidth}
\centering
\includegraphics[height=2in]{vww-v.png}
	\begin{tikzpicture}
	\filldraw 
	(16,4) circle (3pt) node[align=center, below=0.1 cm]{$v_-$}
	(16,6) circle (3pt) node[align=center, below=0.1 cm]{$v_+$}
	(12,4) circle (1pt) node[align=center, below=0.1 cm]{$v_+ w_-$}
	(13,2) circle (1pt) node[align=center, below=0.1 cm]{$v_-  w_-$}
	(13,6) circle (1pt) node[align=center, below=0.1 cm]{$v_+  w_+$}
	(14,4) circle (1pt) node[align=center, below=0.1 cm]{$v_- w_+$}
	(6,2) circle (3pt) node[align=center, below=0.1 cm]{$v_+ w_-  w_-$}
	(6,4) circle (1pt) node[align=center, below=0.1 cm]{$v_+  w_+  w_-$}
	(8,0) circle (3pt) node[align=center, below=0.1 cm]{$v_-  w_-  w_-$}
	(8,2) circle (1pt) node[align=center, below=0.1 cm]{$v_-  w_+  w_-$}
	(8,4) circle (1pt) node[align=center, below=0.1 cm]{$v_+  w_-  w_+$}
	(8,6) circle (3pt) node[align=center, below=0.1 cm]{$v_+  w_+  w_+$}
	(10,2) circle (1pt) node[align=center, below=0.1 cm]{$v_- w_-  w_+$}
	(10,4) circle (3pt) node[align=center, below=0.1 cm]{$v_-  w_+  w_+$};
	\path[->] (8.2,6) edge (12.8,6);
	\path[->] (13.2,6) edge (15.8,6);
	\path[->] (10.2,4) edge[bend left] (13.8,4);
	\path[->] (14.2,4) edge[bend left] (15.8,4);
	\path[->] (6.2,2) edge[bend left] (12.8,2);
	\end{tikzpicture}
	\caption{Type $\mathscr{D}$: $ \tilde v \tilde w \to \tilde v$ for Khovanov differentials.}
	\label{kh-vww-v}
\end{minipage} 
\efig

\bfig
\begin{minipage}{\textwidth}
\centering
\includegraphics[height=2in]{ww-vv.png}
	\begin{tikzpicture}
	\filldraw 
	(0,0) circle (3pt) node[align=center, below=0.1 cm]{$w_-w_-$}
	(-1,2) circle (1pt) node[align=center, below=0.1 cm]{$w_+w_-$}
	(1,2) circle (1pt) node[align=center, below=0.1 cm]{$w_-w_+$}
	(0,4) circle (3pt) node[align=center, below=0.1 cm]{$w_+w_+$}
	(3,2) circle (1pt) node[align=center, below=0.1 cm]{$w_-$}
	(3,4) circle (1pt) node[align=center, below=0.1 cm]{$w_+$}
	(5,4) circle (3pt) node[align=center, below=0.1 cm]{$v_+v_-$}
	(6,2) circle (3pt) node[align=center, below=0.1 cm]{$v_-v_-$}
	(6,6) circle (3pt) node[align=center, below=0.1 cm]{$v_+v_+$}
	(7,4) circle (3pt) node[align=center, below=0.1 cm]{$v_-v_+$};
	\path[->] (0.2,4) edge (2.8,4);
	\path[->] (3.2,4) edge[bend left] (4.8,4);
	\path[->] (3.2,4) edge[bend left] (6.8,4);
	\end{tikzpicture}
	\caption{Type $\mathscr{B}$: $\tilde w \to \tilde v \tilde v$ for Khovanov differentials.}
	\label{kh-ww-vv}
\end{minipage}

\begin{minipage}{\textwidth}
\centering
\hspace{-1in}\includegraphics[height=2in]{vv-ww.png}
	\begin{tikzpicture}
	\filldraw 
	(7,4) circle (3pt) node[align=center, below=0.1 cm]{$w_-w_-$}
	(6,6) circle (1pt) node[align=center, below=0.1 cm]{$w_+w_-$}
	(8,6) circle (1pt) node[align=center, below=0.1 cm]{$w_-w_+$}
	(7,8) circle (3pt) node[align=center, below=0.1 cm]{$w_+w_+$}
	(4,4) circle (1pt) node[align=center, below=0.1 cm]{$w_-$}
	(4,6) circle (1pt) node[align=center, below=0.1 cm]{$w_+$}
	(0,4) circle (3pt) node[align=center, below=0.1 cm]{$v_+v_-$}
	(1,2) circle (3pt) node[align=center, below=0.1 cm]{$v_-v_-$}
	(1,6) circle (3pt) node[align=center, below=0.1 cm]{$v_+v_+$}
	(2,4) circle (3pt) node[align=center, below=0.1 cm]{$v_-v_+$};
	\path[->] (0.2,4) edge[bend left] (3.8,4);
	\path[->] (4.2,4) edge (6.8,4);
	\path[->] (2.2,4) edge[bend left] (3.8,4);
	\path[->] (1.2,6) edge (3.8,6);
	\path[->] (4.2,6) edge[bend left] (5.8,6);
	\path[->] (4.2,6) edge[bend left] (7.8,6);
	\end{tikzpicture}
	\caption{Type $\mathscr{E}$: $\tilde v \tilde v \to \tilde w$ for Khovanov differentials.}
	\label{kh-vv-ww}
\end{minipage} 
\efig

\bfig
\begin{minipage}{\textwidth}
\centering

\hspace{-1in} \includegraphics[height=2in]{ww-wwww.png}

\vspace{0.5in}

\hspace{-1in}
	\begin{tikzpicture}
	\filldraw 
	(0,2) circle (1pt) node[align=center, below=0.1 cm]{$w_+ w_-$}
	(1,0) circle (3pt) node[align=center, below=0.1 cm]{$w_- w_-$}
	(2,2) circle (1pt) node[align=center, below=0.1 cm]{$w_-w_+$}
	(1,4) circle (3pt) node[align=center, below=0.1 cm]{$w_+w_+$}
	(4,2) circle (1pt) node[align=center, below=0.1 cm]{$w_+w_-w_-$}
	(4,4) circle (1pt) node[align=center, below=0.1 cm]{$w_+w_+w_-$}
	(6,0) circle (1pt) node[align=center, below=0.1 cm]{$w_-w_-w_-$}
	(6,2) circle (1pt) node[align=center, below=0.1 cm]{$w_-w_+w_-$}
	(6,4) circle (1pt) node[align=center, below=0.1 cm]{$w_+w_-w_+$}
	(6,6) circle (1pt) node[align=center, below=0.1 cm]{$w_+w_+w_+$}
	(8,2) circle (1pt) node[align=center, below=0.1 cm]{$w_-w_-w_+$}
	(8,4) circle (1pt) node[align=center, below=0.1 cm]{$w_-w_+w_+$}
	(10,4)  circle (3pt) node[align=center, below=0.1 cm]{$w_+w_+ w_-w_-$}
	(12.5,0) circle (3pt) node[align=center, below=0.1 cm]{$w_- w_-w_- w_-$}
	(12.5,8) circle (3pt) node[align=center, below=0.1 cm]{$w_+w_+w_+w_+$}
	(15,4) circle (3pt) node[align=center, below=0.1 cm]{$w_- w_-w_+w_+$}
	(11,2) circle (1pt) node{}
	(11,4) circle (1pt) node{}
	(11,6) circle (1pt) node{}
	(12,2) circle (1pt) node{}
	(12,4) circle (1pt) node{}
	(12,6) circle (1pt) node{}
	(13,2) circle (1pt) node{}
	(13,4) circle (1pt) node{}
	(13,6) circle (1pt) node{}
	(14,2) circle (1pt) node{}
	(14,4) circle (1pt) node{}
	(14,6) circle (1pt) node{};
	\path[->] (1.2,4) edge[bend left] (4-.2,4);
	\path[->] (1.2,4) edge[bend left] (8-.2,4);
	\path[->] (4.2,4) edge[bend left] (10-.2,4);
	\path[->] (4.2,4) edge[bend left] (12-.2,4);
	\path[->] (8.2,4) edge[bend left] (13-.2,4);
	\path[->] (8.2,4) edge[bend left] (15-.2,4);
	\path[->] (1.2,0) edge (5.8,0);
	\path[->] (6.2,0) edge (12.5-.2,0);
	\end{tikzpicture}
	\caption{Type $\mathscr{C}$: $\tilde w \to \tilde w \tilde w$ for Khovanov differentials.}
	\label{kh-ww-wwww}
\end{minipage} 
\efig

\bfig
\begin{minipage}{\textwidth}
\centering

\includegraphics[height=2in]{wwww-ww.png}

\vspace{0.5in}

\hspace{-1in}
	\begin{tikzpicture}
	\filldraw 
	(0+22,2+4) circle (1pt) node[align=center, below=0.1 cm]{$w_+ w_-$}
	(1+22,0+4) circle (3pt) node[align=center, below=0.1 cm]{$w_- w_-$}
	(2+22,2+4) circle (1pt) node[align=center, below=0.1 cm]{$w_-w_+$}
	(1+22,4+4) circle (3pt) node[align=center, below=0.1 cm]{$w_+w_+$}
	(4+12,2+2) circle (1pt) node[align=center, below=0.1 cm]{$w_+w_-w_-$}
	(4+12,4+2) circle (1pt) node[align=center, below=0.1 cm]{$w_+w_+w_-$}
	(6+12,0+2) circle (1pt) node[align=center, below=0.1 cm]{$w_-w_-w_-$}
	(6+12,2+2) circle (1pt) node[align=center, below=0.1 cm]{$w_-w_+w_-$}
	(6+12,4+2) circle (1pt) node[align=center, below=0.1 cm]{$w_+w_-w_+$}
	(6+12,6+2) circle (1pt) node[align=center, below=0.1 cm]{$w_+w_+w_+$}
	(8+12,2+2) circle (1pt) node[align=center, below=0.1 cm]{$w_-w_-w_+$}
	(8+12,4+2) circle (1pt) node[align=center, below=0.1 cm]{$w_-w_+w_+$}
	(10-1,4)  circle (3pt) node[align=center, below=0.1 cm]{$w_+w_+ w_-w_-$}
	(12.5-1,0) circle (3pt) node[align=center, below=0.1 cm]{$w_- w_-w_- w_-$}
	(12.5-1,8) circle (3pt) node[align=center, below=0.1 cm]{$w_+w_+w_+w_+$}
	(15-1,4) circle (3pt) node[align=center, below=0.1 cm]{$w_- w_-w_+w_+$}
	(11-1,2) circle (1pt) node{}
	(11-1,4) circle (1pt) node{}
	(11-1,6) circle (1pt) node{}
	(12-1,2) circle (1pt) node{}
	(12-1,4) circle (1pt) node{}
	(12-1,6) circle (1pt) node{}
	(13-1,2) circle (1pt) node{}
	(13-1,4) circle (1pt) node{}
	(13-1,6) circle (1pt) node{}
	(14-1,2) circle (1pt) node{}
	(14-1,4) circle (1pt) node{}
	(14-1,6) circle (1pt) node{};
	\path[->] (12.5-1+0.2,8) edge (6+12-0.2,6+2);
	\path[->] (6+12+0.2,6+2) edge (1+22-0.2,4+4);
	\path[->] (10-1+0.2,4) edge[bend left] (6+12-0.2,2+2);
	\path[->] (6+12+0.2,2+2) edge[bend left]  (0+23-0.2,4);
	\path[->] (15-1+0.2,4) edge[bend left]  (8+12-0.2,2+2);
	\path[->] (8+12+0.2,2+2) edge[bend left]  (0+23-0.2,4);

	\end{tikzpicture}
	\caption{Type F: $\tilde w \tilde w \to \tilde w$ for Khovanov differentials. }
	\label{kh-wwww-ww}
\end{minipage}

\efig
\epf

\eprop

Since the proof that higher differentials vanish relies on the $k$-grading, we cannot use the same methods to determine whether higher differentials vanish in the $\khtate$ case. However, the $k$-grading is only needed when considering atoms of type $\SE$ (see \S \ref{grading-shifts-of-higher-differentials}), so one can still consider special cases where this obstruction is not needed for the conjecture to hold.

\subsection{Positive and negative braid closures}

Recall that the $k$-grading is used in the proof of Theorem \ref{main-thm} only when the $j$-grading was insufficient, namely in the situation where all the first equivariant steps through the subcube are $\SE$ atoms. From this we deduce that Conjecture \ref{main-conj} holds for ``mostly negative'' links.

\begin{reptheorem}{mostly-negative-thm}
If the annular braid closure $L = \widehat \beta$ has a diagram with at most 1 positive crossing, then the spectral sequence in Conjecture \ref{main-conj} exists and the cascade of rank inequalities holds.

\bpf
The proof of Theorem \ref{main-thm} relies mostly on the $j$-grading to prove that the $\hvE$ spectral sequence collapses on page $\hvE^3$, except for the case of the generator $v_+v_+$ in the starting resolution of Type $\SE$ atoms. In the $\khtate$ case, again the $j$-grading is sufficient for all cases except for Type $\SE$, so for cases where at most one atom of Type $\SE$ appears, longer differentials cannot exist.  This corresponds to the case where $Cube(\CD(L))$ has at most one edge corresponding to a change of resolution from two nontrivial circles to two trivial circles. This means that the all-zeros resolution of $\CD(L)$ corresponds to at most one braid-like resolution; that is, the braid $\beta$ has at most one positive crossing.
\epf
\end{reptheorem}

On the other hand, we cannot use the same method to prove the existence of the spectral sequence for ``mostly positive'' braid closures. However, we can use duality in Khovanov and annular Khovanov homology to show that the cascade of rank inequalities still holds in this case.

\begin{repcorollary}{mostly-positive-cor}
If the annular braid closure $L = \widehat \beta$ has a diagram with at most 1 negative crossing, then the cascade of rank inequalities holds. 

\bpf
Let $L$ be a link in $S^3$ and $\CD$ a diagram of $L$. Let $L^!$ denote the mirror of $L$, and $\CD^!$ the dual diagram to $\CD$ (switch the sign of all crossings). Observe that the complex $(CKh(\CD), \partial^{Kh})$ is dual to $(CKh(\CD^!), \partial^{Kh})$. The cohomology of $CKh(\CD)$  is isomorphic to the homology of $CKh(\CD^!)$, which in turn is isomorphic to the cohomology of $CKh(\CD^!)$, as we are working over field coefficients. (It is now important to note that what we have been calling ``Khovanov homology'' is actually a cohomology theory.) Therefore $rk_{\FF} Kh(\CD) = rk_{\FF} Kh(\CD^!)$. For an annular link $L \subset A \times I$, the above paragraph holds analogously for annular Khovanov (co)homology. Therefore it follows from Theorem \ref{mostly-negative-thm} that the cascade of rank inequalities also holds for mostly positive braid closures.
\epf
\end{repcorollary}

\brmk
Aside from proving the conjecture, many related questions remain. Are there other families of links for which the conjecture holds? Do certain steps in the proof for the $\atate$ spectral sequence point to other obstructions to longer differentials for the $\khtate$ spectral sequence? What are some other algebraic conditions on a braid closure which guarantee that the conjecture holds? 
\ermk

\brmk
While Khovanov-thin links have been studied for some time, at present I am unaware of any explorations of ``annular Khovanov-thin links.'' One could define an annular link $L$ to be \emph{$AKh$-thin} if for each pair of quantum and $\sltwo$ weight space gradings $(j,k)$, $AKh^{j,k}$ is supported on two adjacent $i$-gradings. Then $AKh$-thinness of the quotient link $L$ suffices for the spectral sequence in Conjecture \ref{main-conj}, since all differentials after page $\hvE^3$ are too long in $i$-degree to be nontrivial. 
\ermk

\bibliographystyle{abbrv}
\bibliography{biblio2}

\end{document}